\documentclass[11pt,reqno]{article}
\usepackage{hyperref}
\usepackage{color} 

\let\eps=\varepsilon

\usepackage{amsmath}
\usepackage{amssymb}
\usepackage{amsthm}
\usepackage{srcltx}
\usepackage{enumerate}
\numberwithin{equation}{section}

\newtheorem{lemma}{Lemma}[section]
\newtheorem{theorem}[lemma]{Theorem}

\newtheorem{proposition}[lemma]{Proposition}

\newtheorem{remark}[lemma]{Remark}
\newtheorem{remarks}[lemma]{Remarks}

\def\Gammato{\to\!\!\!\!\!\!\!{}^\Gamma\;\;}   
\def\Gammatopi{\to\!\!\!\!\!\!\!\!\!{}^{\Gamma\widetilde\pi}\;\;}   
  
\def\2to{\to\!\!\!\!\!\!{}_{_2}\;\;}
\def\pito{\;\!\!\to\!\!\!\!\!\!\!_{_{\scriptstyle\widetilde\pi}}\;\;}

\def\w2to{\rightharpoonup\!\!\!\!\!\!{}_{_2}\;\;}
\def\ws2to{\rightharpoonup\!\!\!\!\!\!{}_{_2}\!\!\!{}^{*}\;\;\;}
\def\wto{\rightharpoonup}
\def\wsto{\rightharpoonup\!\!\!\!\!\!{}^{*}\;\,}

\def\erre
{{\bf R}}

\def\enne
{{\bf N}}

\def\Theta{{\mit \Theta}}
\let\La=\Lambda\def\Lambda{{\mit \La}}
\let\Sg=\Sigma\def\Sigma{{\mit \Sg}}
\let\Pig=\Pi\def\Pi{{\mit \Pig}}
\let\Ps=\psi\def\psi{{\mit \Ps}}
\let\Xii=\Xi\def\Xi{{\mit \Xii}}
\def\graph{\mathop{\rm graph}}
\def\div{\mathop{\rm div}} 
\def\curl{\mathop{\rm curl}} 
\def\gr{\mathop{\rm graph}} 
\def\lsc{\mathop{\rm sc_w^-}}

\title{\bf Structural Compactness and Stability of Pseudo-Monotone Flows}
\author{Augusto Visintin
\thanks{Dipartimento di Matematica dell'Universit\`a degli Studi di Trento --
via Sommarive 14,  38050 Povo di Trento, Italia -- email: augusto.visintin@unitn.it }
}

\date{\today} 
\begin{document}
\maketitle

\begin{abstract} 
Fitzpatrick's variational representation of maximal monotone operators is here extended 
to a class of pseudo-monotone operators in Banach spaces. 
On this basis, the initial-value problem associated with the first-order flow 
of such an operator is here reformulated as a minimization principle, 
extending a method that was pioneered by Brezis, Ekeland and Nayroles for gradient flows.
This formulation is used to prove that the problem is stable w.r.t.\ arbitrary perturbations 
not only of data but also of operators.
This is achieved by using the notion of {\it evolutionary $\Gamma$-convergence\/} w.r.t.\ a nonlinear topology of weak type. 

These results are applied to the Cauchy problem for quasilinear parabolic PDEs.
This provides the structural compactness and stability of the model of 
several physical phenomena:
nonlinear diffusion, 
incompressible viscous flow,
phase transitions, and so on.
\end{abstract} 

\bigskip 
\noindent{\bf Keywords:} 
Fitzpatrick theory, 
Pseudo-mono\-tone operators,
%Structural compactness,
%Structural stability,  
Variational formulation, 
Nonlinear weak-type topology,
Quasilinear parabolic PDEs. 
 
\bigskip 
\noindent{\bf AMS Classification (2000):}  
35K60, % Boundary value problems for nonlinear parabolic PDE  
47H05, % Monotone operators (with respect to duality)  
49J40, % Variational methods including variational inequalities      
58E. % Variational problems in infinite-dimensional spaces 

\section{Introduction}\label{intro} 

\noindent  
This work extends to a class of nonmonotone operators the variational formulation of maximal monotone operators of \cite{Fi}.
On that basis, first-order flows associated to those operators are here formulated as 
minimization principles.
Their stability w.r.t.\ arbitrary perturbations of data and operators is then proved by using 
a nonlinear topology, and the novel notion of {\it evolutionary $\Gamma$-convergence\/} 
of weak type.

\bigskip
\noindent{\bf Variational representation of maximal monotone flows.}
Let us first outline the variational formulation of first-order flows.
Let $V$ be a reflexive real Banach space, $H$ be a real Hilbert space, and
$V\subset H=H'\subset V'$ with dense injections. 
Let $\varphi: V\to \erre\cup \{+\infty\}$ 
be a lower semi-continuous convex functional, and consider the gradient flow  
\begin{equation}\label{eq.intro.1}
D_tu + \partial\varphi(u)\ni h 
\qquad\text{ in $V'$, a.e.\ in time }(D_t := {\partial/\partial t}).
\end{equation}
By using the classical {\it Fenchel inequality,\/} see e.g.\ \cite{EkTe} or \cite{Fe}, 
Brezis and Ekeland \cite{BrEk} and Nayroles \cite{Na} 
(who actually assumed $V$ to be a Hilbert space)
reformulated this inclusion as the minimization of the nonnegative functional
\begin{equation} \label{eq.intro.2}
\Phi(v,v^*) = \int_0^T [\varphi(v) + \varphi^*(v^* -D_tv) - \langle v^*,v\rangle] \, dt 
+ {1\over2} \|v(T)\|_H^2 - {1\over2} \|u(0)\|_H^2, 
\end{equation} 
under the hypothesis that the minimal value of $\Phi$ is null.
This may be regarded as a variational representation of the monotone operator 
$D_t + \partial\varphi$.

Several years after the works \cite{BrEk} and \cite{Na}, in \cite{Fi} Fitzpatrick introduced 
a {\it representation\/} of maximal monotone operators $\alpha: V\to {\cal P}(V')$ via a minimization principle; this is here briefly reviewed in Section~\ref{sec.fitzp}.
If $f_\alpha: V \!\times\! V' \to \erre\cup \{+\infty\}$ represents the operator $\alpha$, 
one can then extend the formulation 
\eqref{eq.intro.2} to the flow 
\begin{equation}\label{eq.intro.3}
D_tu + \alpha(u)\ni h 
\qquad\text{ in $V'$, a.e.\ in time, }
\end{equation} 
by minimizing the functional
\begin{equation} \label{eq.intro.3-}
\widetilde\Phi(v,v^*) = \int_0^T [f_\alpha(v,v^* -D_tv)- \langle v^*,v\rangle] \, dt 
+ {1\over2} \|v(T)\|_H^2 - {1\over2} \|u(0)\|_H^2. 
\end{equation}
 
\medskip
\noindent{\bf Variational representation of semi-monotone flows.}
The aim of the present paper is twofold.
First, we extend the Fitzpatrick theory to {\it semi-monotone\/} operators 
$\alpha: V\to {\cal P}(V')$. 
These are pseudo-monotone operators of the form $\alpha(v) = \beta(v,v)$ for any $v\in V$,
where $\beta: V^2\to {\cal P}(V')$ is weakly lower semi-continuous w.r.t.\ the first argument 
and maximal monotone w.r.t.\ the second one.   
Examples are often encountered in the theory of nonlinear quasilinear elliptic P.D.E.s.

We represent semi-monotone operator $\alpha$ by means of an (in general nonconvex)
lower semicontinuous function $f:V \!\times\! V'\to \erre\cup \{+\infty\}$.
The lack of convexity induces us to deal with a special nonlinear topology of weak-type.
We then formulate the corresponding flow $D_tu + \alpha(u)\ni h$ as in \eqref{eq.intro.3-}. 
Existence of a solution for the associated Cauchy problem can be proved via the 
variational formulation.
However, structural properties are the main concern of this work. 

\bigskip
\noindent{\bf Structural compactness and structural stability.}
Aside the customary notion of {\it stability,\/} i.e.\ robustness to perturbations of the data
(e.g., source terms, initial- and boundary-values),
one may consider {\it structural stability,\/} 
i.e., robustness to perturbations of the structure of the problem
(e.g.\ linear and nonlinear operators in differential equations).
These notions have an obvious applicative motivation: 
data and operators are accessible just with some approximation;
moreover perturbations, e.g.\ in the coefficients of quasilinear operators, 
are hardly predictable. 
If properly formulated, structural stability may thus be regarded as a basic requisite   
for the feasibility, or even the applicative soundness, of mathematical models. 
Structural stability is also a prerequisite for numerical treatment and efficient control
of systems. 

Structural stability is here associated with {\it structural compactness,\/}
namely the existence of a convergent subsequence of the solutions of the perturbed problems.  
In order to fix ideas, let us consider a model problem of the form $Au\ni h$, 
$A$ being a multi-valued operator acting between Banach spaces and $h$ a datum.
Given bounded families of data and of operators, we formulate the stability of the problem
in terms of two properties:

(i) {\it structural compactness:\/} existence of a pair $(A,h)$ and of sequences 
of data $\{h_n\}$ and of operators $\{A_n\}$ such that
$A_n\to A$ and $h_n\to h$ w.r.t.\ topologies that must be specified;

(ii) {\it structural stability:\/} if $A_nu_n\ni h_n$ for any $n$, 
$A_n\to A$, $h_n\to h$ and $u_n\to u$, 
then $u$ is a solution of the asymptotic problem $Au\ni h$.  

For variational principles, De Giorgi's theory of $\Gamma$-convergence fulfills both 
properties, see \cite{DeFr} and e.g.\ the monographs 
\cite{At},\cite{Bra1},\cite{Bra2},\cite{Da}. 
This typically applies to stationary models.
The present article extends that approach to flows for semi-monotone operators,
by using the above variational formulation and the novel notion of 
{\it evolutionary $\Gamma$-convergence\/} of weak type, see \cite{Vi17}. 
 
Several authors have been dealing with structural stability, but not with structural compactness,
as far as this author knows.
In several cases, see e.g.\ \cite{Mi1},\cite{RoRoSt},\cite{St1}, 
stability has been addressed via so-called Mosco-convergence, 
i.e.\ $\Gamma$-convergence with respect to both the weak and the strong topology, 
see e.g.\ \cite{At},\cite{Mo}.
One might wonder whether that notion is feasible in the present set-up, too.
Although the Mosco-convergence has a wealth of properties that make it quite attractive 
for the analyst, this author does not see how it might to applied when dealing 
with structural compactness, since results of Mosco-compactness seem to be quite rare.
On the other hand $\Gamma$-compactness holds under minimal restrictions.  

\bigskip
\noindent{\bf Nonlinear weak topology and evolutionary $\Gamma$-convergence.} 
As we pointed out above, the choice of the topology for data and operators 
plays a key role in the analysis of structural properties:
not only the $\Gamma$-convergence but also the $\Gamma$-limit typically 
obviously depend on the selected topology.
The point that we advocate in this work is that the relevant topology should not be 
left to the convenience of mathematical analysis 
(e.g., the exigence of passing to the limit in nonlinear terms), 
but should depend on a selection criterion rooted in the model.
As a criterion we suggest the occurrence of structural compactness.
By this we mean that the topology must be so weak that,
under minimal assumptions of either bounded or coerciveness, 
there exists a sequence of solutions of the perturbed problems that converges
with respect to that topology. 
In order to fulfill a compactness theorem, this topology is necessarily of weak type.
However the weak topology is not appropriate, since it would not provide the 
existence of a recovery sequences: a slightly stronger topology is needed.
The solution is provided by what we name the 
{\it nonlinear weak topology\/} of $V \!\times\! V'$, see \cite{ViCalVar}. 
\footnote{ Here we are concerned with pseudo-monotone operators, 
for which rather weak a priori estimates are just available.
In case of monotone operators, stronger estimates may be proved 
under further assumptions on the data.
The monotonicity assumption allows one to deal with structural properties 
in smaller spaces w.r.t.\ a stronger topology, 
without the need of using our nonlinear weak topology. 
The corresponding analysis would then be simpler than the present one.
}

Here we show the structural stability of the variational formulation of of flows 
of semi-monotone operators.
The analogous property fails for the class of gradient flows (like \eqref{eq.intro.1}): 
in Section~5 of \cite{ViCalVar} it was shown that the class of 
subdifferentials of lower semicontinuous convex functions is not stable by 
$\Gamma$-convergence w.r.t.\ the nonlinear weak topology. 

Our analysis is based on an extension of the notion of $\Gamma$-convergence
to evolutionary problems; this is performed in the parallel work \cite{Vi17},
and is here illustrated in Section~\ref{sec.evol}.
This definition is not equivalent either to that of \cite{SaSe} or to that of 
\cite{DaSa}, \cite{Mi1}, \cite{Mi2}.
In those works $\Gamma$-convergence is actually assumed for almost any 
$t\in {}]0,T[$, whereas here it is just weak in $L^1_\mu(0,T)$. 

The present definition fits the rather general framework of $\bar\Gamma$-convergence, 
that is defined in Chap.~16 of \cite{Da}, see also references therein.
The results of that monograph however do not encompass the theorem of
$\Gamma$-compactness of Section~\ref{sec.evol}.  

\bigskip
\noindent{\bf Plan of work.} 
In Section~\ref{sec.fitzp} we review the tenets of the Fitzpatrick theory;
as an example, in Proposition~\ref{eq.fitp.ell} we {\it represent\/} a quasilinear elliptic 
differential operator.
In Section~\ref{sec.repres} we extend that theory to nonmonotone operators, and
represent sums of operators in Theorems~\ref{sumcor} and \ref{sum}. 
In Section~\ref{sec.BEN} we introduce the class of semi-monotone operators, 
and represent them in Theorem~\ref{BEN.param} and Proposition~\ref{eq.fitp.ellpse}.
In Section~\ref{sec.time} we extend the previous results to spaces of time-dependent functions,
see Theorem~\ref{timeint}.
On that basis in Section~\ref{sec.psflow} we generalize the 
{\it Brezis-Ekeland-Nayroles principle\/} (``BEN principle'' for short)
to the first-order flow of semi-monotone operators;
see Theorems~\ref{prop.repr} and \ref{prop.repr'}.

In Section~\ref{sec.comstab} we review the nonlinear weak topology $\widetilde\pi$, and
illustrate the structural compactness and stability of the variational formulation 
of first-order flows, see Theorems~\ref{comsta} and \ref{teo.comp'}.  
In Section~\ref{sec.par} we then apply the two latter theorems to 
quasilinear equations: in Theorem~\ref{eq.par.stru} we show the structural compactness 
and structural stability of the flow \eqref{eq.intro.3}.
In Section~\ref{sec.PDE} we illustrate a number of quasilinear PDEs of mathematical physics, which the previous results apply to. 

Finally, in the Appendix (Section~\ref{sec.evol}) we outline the notion of 
evolutionary $\Gamma$-convergence of weak type,
and state a result of \cite{Vi17} which is used in the proof of Theorem~\ref{teo.comp'}. 

\bigskip
\noindent{\bf A look at the literature.} 
The note \cite{Kr} of Krylov of 1982 perhaps was the first to address
the variational formulation of maximal monotone operators (in finite dimensional spaces). 
Fitzpatrick's Theorem~\ref{teo.Fi} appeared in 1988 in \cite{Fi},
and was not noticed for several years. 
It was eventually rediscovered by Martinez-Legaz and Th\'era \cite{MaTh01}
and (independently) by Burachik and Svaiter \cite{BuSv02} in the early 2000s. 
Since then this theory has been continuously progressing, see e.g.\
\cite{BaBoWa}, \cite{BaWa}, \cite{Bor}, \cite{BuSv03},  
\cite{MaSv05}, \cite{MaSv08}, \cite{Pe04}, \cite{Pe04'}, just to mention few contributions. 
This author addressed this theory for monotone operators e.g.\ in 
\cite{ViAMSA}--\cite{ViCalVar}. 
The representation of (single-valued) nonmonotone operators was first outlined in \cite{Vi14}. 

Fitzpatrick's formulation of maximal monotone operators via a minimization principle cast 
a new light upon the formulation of gradient flows that Brezis and Ekeland \cite {BrEk} and 
Nayroles \cite{Na} had pointed out in 1976, prior to the Fitzpatrick theorem 
of 1988. The BEN principle was studied and applied in several papers; 
see for instance \cite{Auc}, \cite{Ri1}, \cite{RoRoSt}, \cite{Rou2} and \cite{St1}.
The extension \eqref{eq.intro.3-} to maximal monotone flows of the form \eqref{eq.intro.3} 
was introduced in \cite{ViAMSA}.   

Single-valued semi-monotone operators were introduced by Browder in \cite{Bro1}.
In \cite{Ke} (see pp.~220-222) Kenmochi extended this notion to 
multi-valued operators, and applied it to elliptic problems.
That definition is here amended in a form that is convenient for the analysis of 
quasilinear parabolic equations.
These operators are comprised in the larger class of (multi-valued)
{\it generalized pseudo-monotone\/} operators, in the sense of Browder-Hess \cite{BrHe}; 
see also e.g.\ \cite[p.\ 368]{HuPa}. This latter class extends Brezis's notion of (single-valued) 
pseudo-monotone operators of \cite{Br0}.
Although semi-monotone operators are less general than those of  \cite{BrHe}, 
they include many examples of applicative interest. 

A comparison of the present results with those of \cite{ViCalVar}, \cite{Vi15} and 
\cite{ViGil} is in order. 
\cite{ViCalVar} dealt with structural compactness and stability of maximal monotone flows,
but excluded the onset of long memory in the $\Gamma$-limit only for a restricted class 
of equations.  
The present work deals with semi-monotone operators and uses a different approach, 
that rests upon the notion of evolutionary $\Gamma$-convergence of weak type.  
\cite{Vi15} addressed the structural stability of \eqref{eq.intro.3}, 
but dealt with single-valued pseudo-monotone operators and used a different approach, 
unrelated to the BEN-Fitzpatrick theory.
\cite{ViGil} announced the present results.

The main novelties of the present work include:
(i) the extension of the Fitzpatrick theory to nonmonotone operators, 
in particular semi-monotone operators;
(ii) the extension of the BEN minimization principle to the corresponding flows;
(iii) the general analysis of structural compactness and structural stability 
via a nonlinear topology of weak type;
(iv) the application of those results to the Cauchy problem for quasilinear parabolic PDEs. 
 
\section{Outline of Fitzpatrick's theory}
\label{sec.fitzp}

\noindent 
In this section we review a pioneering result of Fitzpatrick \cite{Fi}  
for monotone operators, and outline the related theory. 
 
\bigskip
\noindent{\bf The Fitzpatrick theorem.} 
Let $V$ be a real Banach space with norm $\|\cdot\|$, dual norm $\|\cdot\|_{V'}$, and
duality pairing $\langle \cdot,\cdot \rangle$. 
Let $\alpha: V \to {\cal P}(V')$ be a (possibly multi-valued) measurable operator,
i.e., such that $g^{-1}(A) := \big\{v\in V: g(v)\cap A\not= \emptyset \big\}$
is measurable, for any open subset $A$ of $V'$. 
For instance, this condition is fulfilled if $\alpha$ is maximal monotone.
We shall always assume that $\alpha$ is proper, i.e., $\alpha(V) \not= \emptyset$. 
In \cite{Fi} Fitzpatrick defined what is now called the {\it Fitzpatrick function:\/}
\begin{equation}\label{eq.fitzp.Fitzfun}
\begin{split}
f_\alpha(v,v^*) :=
& \; \langle v^*,v\rangle + \sup \big\{ \langle v^*- \widetilde v^*,\widetilde v - v\rangle:
\widetilde v\in V, \widetilde v^*\in \alpha(\widetilde v) \big\} 
\\[2mm]
=& \sup \big\{ \langle v^*, \widetilde v\rangle - \langle \widetilde v^*, \widetilde v- v\rangle:
\widetilde v\in V,  \widetilde v^*\in \alpha(\widetilde v) \big\} 
\quad\forall (v,v^*)\in V \!\times\! V'.
\end{split}  
\end{equation} 
This function is convex and lower semi-continuous, since it is
the supremum of a family of affine and continuous functions.

\begin{theorem}[\cite{Fi}]\label{teo.Fi}
An operator $\alpha: V \to {\cal P}(V')$ is maximal monotone if and only if
\begin{eqnarray}
&f_\alpha(v,v^*) \ge \langle v^*,v\rangle
\qquad\forall (v,v^*)\in V \!\times\! V', 
\label{eq.fitzp.Fitztheo1}
\\[2mm]
&f_\alpha(v,v^*) = \langle v^*,v\rangle
\quad\Leftrightarrow\quad v^*\in \alpha(v). 
\label{eq.fitzp.Fitztheo2}
\end{eqnarray} 
\end{theorem}

This theorem generalizes the following classical result of Fenchel. 
Let $\varphi: V\to \erre\cup \{+\infty\}$ be a 
(possibly nonconvex and non-lower-semi-continuous) proper function
(i.e., $\varphi\not\equiv +\infty$), and denote by
$\varphi^*: V'\to \erre\cup \{+\infty\}$ and $\partial\varphi: V\to {\cal P}(V')$ 
respectively the convex conjugate function and the subdifferential of $\varphi$;
see e.g.\ \cite{EkTe}, \cite{Fe}. 
Then
\begin{eqnarray}
&\varphi(v)+ \varphi^*(v^*) \ge \langle v^*,v\rangle
\qquad\forall (v,v^*)\in V \!\times\! V', 
\label{eq.fitzp.Fen1}
\\[2mm]
&\varphi(v)+\varphi^*(v^*)= \langle v^*,v\rangle
\;\;\Leftrightarrow\;\;
v^*\in\partial \varphi(v). 
\label{eq.fitzp.Fen2}
\end{eqnarray}
We shall refer to \eqref{eq.fitzp.Fitztheo1} and \eqref{eq.fitzp.Fitztheo2} 
(\eqref{eq.fitzp.Fen1} and \eqref{eq.fitzp.Fen2}, resp.) as the 
{\it Fitzpatrick system\/} (the {\it Fenchel system,\/} resp.), and to the mapping
$(v,v^*)\mapsto \varphi(v)+\varphi^*(v^*)$ as the {\it Fenchel function\/}
of the operator $\partial\varphi$.
 
\bigskip 
\noindent{\bf Representative functions.}
Extending the theory of Fitzpatrick, nowadays one says that a function $f$ (variationally) 
{\it represents\/} a proper measurable operator $\alpha:V \to {\cal P}(V')$ whenever
\begin{eqnarray}  
&f:V \!\times\! V' \to \erre\cup \{+\infty\} 
\text{ \ is convex and lower semi-continuous, }
\label{eq.fitzp.convrep1}
\\
&f(v,v^*) \ge \langle v^*,v\rangle
\qquad\forall (v,v^*)\in V \!\times\! V',
\label{eq.fitzp.convrep2}
\\
&f(v,v^*) = \langle v^*,v\rangle
\quad\Leftrightarrow\quad v^*\in \alpha(v).
\label{eq.fitzp.convrep3} 
\end{eqnarray}

One accordingly says that $\alpha$ is {\it representable,\/}  
that $f$ is a (convex) {\it representative function,\/} 
that $f$ {\it represents\/} $\alpha$, and so on.
For instance, $\partial \varphi$ represents $\varphi$,
because of \eqref{eq.fitzp.Fen1} and \eqref{eq.fitzp.Fen2}. 
We shall denote by ${\cal F}(V)$ the class of the functions that 
fulfill \eqref{eq.fitzp.convrep1} and \eqref{eq.fitzp.convrep2}.
Representable operators are monotone, see \cite{Fi}, 
but need not be either cyclically monotone or maximal monotone.
However, by Theorem~\ref{teo.Fi} any operator that is represented by its 
Fitzpatrick function is maximal monotone. 

\bigskip
\noindent{\bf Examples.}
(i) For any (possibly nonconvex and non-lower-semi-continuous) proper function 
$\varphi: V\to \erre\cup \{+\infty\}$, the Fenchel function
$f: (v,v^*)\mapsto \varphi(v)+ \varphi^*(v^*)$ represents the operator $\partial\varphi$. 

\smallskip
(ii) Let $A:V\to V'$ be a linear, bounded and invertible monotone operator, 
and define the convex and continuous mapping
\begin{equation}\label{eq.fitzp.counterex}
F_b: V\times V'\to \erre: (v,v^*)\mapsto b[\langle Av,v\rangle + \langle v^*,A^{-1}v^*\rangle]
\qquad\forall b> 0.
\end{equation}
$F_{1/2}$ is the Fenchel function of the operator $A$.
For any $b> 1/2$, $F_b$ represents the monotone but nonmaximal monotone operator 
$\alpha(0) =\{0\}$, $\alpha(v) = \emptyset$ for any $v\not=0$.
For $0< b< 1/2$, $F_b$ does not represents any operator, as \eqref{eq.fitzp.convrep2} fails. 

\smallskip
(iii) First let us denote by $I_\alpha$ the indicator function
of the graph of any operator $\alpha:V \to {\cal P}(V')$:
\[
I_\alpha(v,v^*) =0 \quad\mbox{ if }v^*\in \alpha(v),
\qquad
I_\alpha(v,v^*) =+\infty \quad\mbox{ if }v^*\not\in \alpha(v).
\]

The Fitzpatrick function of the operator $A$ of example (ii) reads
\begin{equation}\label{eq.fitzp.counterex=}
f_A(v,v^*) = I_A(v,v^*) + \langle Av,v\rangle \;
(= I_A(v,v^*) + \langle v^*,v\rangle)
\qquad\forall (v,v^*)\in V \!\times\! V'. 
\end{equation}
This function is not easily handled. It seems more convenient to represent linear maximal monotone operators by the Fenchel function $F_{1/2}$, that we just defined.
   
\smallskip
(iv) This example is of paramount importance for applications to PDEs.
Let $\Omega$ be a bounded domain of $\erre^N$ ($N\ge1$), $p\in{} [2,+\infty[$, 
and set $p'=p/(p-1)$ and 
\begin{equation}\label{eq.fitzp.spaces}
V:= W^{1,p}_0(\Omega), \quad
H= L^2(\Omega) =H', \qquad
\text{ whence }\qquad V' = W^{-1,p'}(\Omega). 
\end{equation}
Let a maximal monotone mapping $\vec\gamma: \erre^N\to {\cal P}(\erre^N)$
be represented by a function $f\in {\cal F}(\erre^N)$, and be such that
\begin{equation}\label{eq.fitp.growth}
\exists c_1,c_2\in\erre^+: \forall \vec w\in\erre^N,
\forall \vec z\in \vec \gamma(\vec w), \qquad
|\vec z| \le c_1 |\vec w|^{p-1} +c_2.
\end{equation} 

This entails that the operator
$\beta: v\mapsto - \nabla \cdot \vec\gamma(\nabla v)$ maps $V$ to ${\cal P}(V')$. 
We claim that $\beta$ is maximal monotone.
Indeed, for any $v^*\in W^{-1,p'}(\Omega)$, the problem
\[
v\in W^{1,p}_0(\Omega),
\qquad
- \nabla \cdot \vec\gamma(\nabla v) - \nabla \cdot (|\nabla v|^{p-2} \nabla v) =v^*
\quad\text{ in }{\cal D}'(\Omega)
\] 
has a solution.
As $v\mapsto - \nabla \cdot (|\nabla v|^{p-2} \nabla v)$
is a duality mapping $W^{1,p}_0(\Omega)\to W^{-1,p'}(\Omega)$,
by the Minty-Browder theorem we conclude that $\beta$ is maximal monotone.

Next we exhibit a representative function of the operator $\beta$. 
Let us first define the linear and continuous operator
$\Lambda: W^{-1,p'}(\Omega)\to W^{1,p'}_0(\Omega)$:
\begin{equation}\label{eq.fitp.du} 
\eta = \Lambda(v^*) 
\quad\Leftrightarrow\quad
\eta \in W^{1,p'}_0(\Omega), \;\;
- \Delta \eta = v^* \text{ in }{\cal D}'(\Omega).
\end{equation}  

\begin{proposition}\label{eq.fitp.ell} 
Let a maximal monotone mapping $\vec\gamma: \erre^N\to {\cal P}(\erre^N)$
fulfill \eqref{eq.fitp.growth} and be represented by a function $g\in {\cal F}(\erre^N)$.
Let $\Omega$ be a bounded domain of $\erre^N$ ($N\ge 1$) of Lipschitz class, 
$p\in{} [2,+\infty[$, and set $V:= W^{1,p}_0(\Omega)$. 
Then the maximal monotone operator 
$\beta: V\to {\cal P}(V'): v\mapsto - \nabla \cdot \vec\gamma(\nabla v)$
is represented by the function 
\begin{equation}\label{eq.fitzp.a2}
\psi(v,v^*) = \int_\Omega g(\nabla v, \nabla \Lambda v^*) \, dx
\qquad\forall (v,v^*)\in V \!\times\! V'.
\end{equation}  
\end{proposition} 

\noindent{\bf Proof.\/}   
The operator $\psi$ is convex, 
because of the convexity of $g$ and of the linearity of the mapping 
$V \!\times\! V'\to L^p(\Omega) \!\times\! L^{p'}(\Omega): 
(v,v^*) \mapsto (\nabla v,\nabla \Lambda v^*)$. 
Notice that 
\begin{equation}\label{eq.fitzp.a3}
\begin{split} 
\varphi(v,v^*) := 
\psi(v,v^*) - \langle v^*, v\rangle 
= \int_\Omega [g(\nabla v,\nabla \Lambda v^*)
- (\nabla\Lambda v^*) \cdot \nabla v]  \, dx &
\\
\forall (v,v^*)\in L^p(\Omega) \!\times\! L^{p'}(\Omega), & 
\end{split} 
\end{equation} 

By \eqref{eq.fitzp.convrep2} the integrand of \eqref{eq.fitzp.a3} is nonnegative, 
hence $\varphi(v,v^*)$ vanishes iff the integrand vanishes a.e.\ in $\Omega$. 
By \eqref{eq.fitzp.convrep3}, this occurs iff 
$\nabla \Lambda v^*\in \vec\gamma(\nabla v)$ a.e.\ in $\Omega$. Hence
\[ 
\psi(v,v^*) = \langle v^*, v\rangle 
\quad\Leftrightarrow\quad
v^* \in - \nabla \!\cdot\! \vec\gamma(\nabla v) = \beta(v) \;\text{ in }{\cal D}'(\Omega).
\] 

By Fatou's lemma and the continuity of the duality pairing, 
$\varphi$ is lower semi-continuous. The same then holds for $\psi$,
which thus represents the operator $\beta$. 
\hfill$\Box$  

\bigskip 
\noindent{\bf Remarks.}
(i) Proposition~\ref{eq.fitp.ell} is easily extended if $\vec\gamma$ (and then $f$) depends 
explicitly on $x$.
In this case the formulation in terms of the single-valued function $f$ avoids technicalities 
related to the measurability of $x$-dependent multi-valued mappings.  

(ii) Further examples of representative functions are provided e.g.\ in \cite{ViCalVar}. 
\hfill$\Box$

\bigskip 
Ahead we shall use the next lemma.

\begin{lemma} [\cite{BuSv02}, \cite{Fi}, \cite{MaSv08}, \cite{Pe04}] \label{lem.dis=}
Let $V$ be a reflexive Banach space, 
$f$ be the Fitzpatrick function of a maximal monotone operator $V\to \erre\cup \{+\infty\}$, 
and $g\in {\cal F}(V)$ be any other representative function of the same operator. 
Then $f\le g$ in $V \!\times\! V'$.
\end{lemma}

\section{Representation of nonmonotone operators}
\label{sec.repres}

\noindent 
After dealing with representative functions of monotone operators,
in this section we extend the representation to general operators, 
and derive some properties. 
The lack of convexity of these functions induces us to deal with the weak topology. 

\bigskip
\noindent{\bf Representation of nonmonotone operators.}   
Let $\alpha: V\to {\cal P}(V')$ be a nontrivial operator 
(i.e., $\alpha(V)\not= \emptyset$).
We shall say that a (possibly nonconvex) function $f$ {\it algebraically represents\/} 
$\alpha$ if 
\begin{eqnarray}  
&f:V \!\times\! V' \to \erre\cup \{+\infty\}, 
\label{eq.fitzp.nonconv1}
\\
&f(v,v^*) \ge \langle v^*,v\rangle
\qquad\forall (v,v^*)\in V \!\times\! V',
\label{eq.fitzp.nonconv2}
\\
&f(v,v^*) = \langle v^*,v\rangle
\quad\Leftrightarrow\quad v^*\in \alpha(v).
\label{eq.fitzp.nonconv3} 
\end{eqnarray} 

Let us set
\begin{equation}\label{eq.pi} 
\pi(v,v^*) = \langle v^*, v\rangle
\qquad\forall (v,v^*)\in V \!\times\! V'. 
\end{equation} 
Obviously, any operator $\alpha: V\to {\cal P}(V')$ is algebraically represented by
$I_\alpha + \pi$. 
This is the largest representative mapping of $\alpha$; 
on the other hand, no smallest representative of $\alpha$ exists.

Next we introduce an algebraical and topological notion of representation.
Because of lack of convexity, we shall distinguish among topologies that 
are intermediate between the weak and the strong topology.
Henceforth we assume that the real Banach space $V$ is reflexive, and say that 
$\rho$ is an {\it intermediate topology\/} on $V \!\times\! V'$ whenever
\begin{equation}\label{fitzp.nonconv4} 
\begin{split} 
&\rho \text{ is stronger than the weak topology of }V \!\times\! V',
\\
&\rho \text{ is weaker than the strong topology of }V \!\times\! V'; 
\end{split}
\end{equation} 
we shall denote this space by $(V \!\times\! V')_\rho$. 
If an operator $\alpha: V\to {\cal P}(V')$ is measurable,
we shall say that $f$ {\it topologically represents\/} $\alpha$ with respect to $\rho$
(or more briefly, $f$ {\it $\rho$-represents\/} $\alpha$) if  
\begin{equation}\label{eq.fitzp.nonconv=}
\text{ $f$ fulfills \eqref{eq.fitzp.nonconv1}--\eqref{eq.fitzp.nonconv3} and is  
lower $\rho$-semi-continuous. }
\end{equation} 

We shall then say that $\alpha$ is  
$\rho$-representable, that it is $\rho$-represented by $f$, and so on. 
Beside the class ${\cal F}(V)$ of convex representatives, see Section~\ref{sec.fitzp},
we shall denote by ${\cal E}(V)$ (${\cal E}_\rho(V)$,  resp.) the class of   
(possibly nonconvex) functions that fulfill 
\eqref{eq.fitzp.nonconv1}--\eqref{eq.fitzp.nonconv3} (\eqref{eq.fitzp.nonconv=}, resp.). 
In particular ${\cal E}_s(V)$ (${\cal E}_w(V)$, resp.) corresponds to $\rho$ equal to the strong (weak, resp.) topology of $V \!\times\! V'$. Thus
\begin{equation}\label{eq.fitzp.inclus}
{\cal F}(V)\subset {\cal E}_w(V) \subset{\cal E}_\rho(V) \subset{\cal E}_s(V)
\subset{\cal E}(V).
\end{equation} 

We shall provide further examples of representative functions ahead. 

\bigskip
\noindent{\bf Some calculus of representative functions.}  
The first part of the next statement will play a role in our extension of the 
Brezis-Ekeland-Nayroles theorem to the flow of nonmonotone operators.

\begin{theorem}\label{sumcor}  
For any operator $\alpha: V\to {\cal P}(V')$ and any positive $L\in {\cal L}(V;V')$,
the following holds:

\indent
(i) If $f$ weakly represents $\alpha$, then $\alpha+L$ is weakly represented by 
\begin{equation}\label{eq.fitzp.sum+}
g_1(v,v^*) = f(v,v^* -Lv) + \langle Lv,v \rangle 
\qquad\forall (v,v^*)\in V \!\times\! V'.
\end{equation} 
\indent
(ii) If $f\in{\cal F}(V)$ then also $g_1\in{\cal F}(V)$.  

\indent
(iii) If the graph of $\alpha$ is weakly closed, 
then the operator $\alpha+L$ is also weakly represented by
\begin{equation}\label{eq.fitzp.sum=}
g_2(v,v^*) = I_\alpha(v,v^* -Lv) + \langle v^*,v \rangle
\qquad\forall (v,v^*)\in V \!\times\! V'. 
\end{equation}
\end{theorem} 

\noindent{\bf Proof.\/}   
(i) As $f$ is weakly lower semi-continuous, $L$ is weakly continuous, 
and the second addendum of \eqref{eq.fitzp.sum+} is convex and weakly continuous, 
it follows that $g_1$ is weakly lower semi-continuous. 
As the properties \eqref{eq.fitzp.nonconv2} and \eqref{eq.fitzp.nonconv3} hold for $f$,
they are also fulfilled by $g_1$.

\smallskip
(ii) If $f\in{\cal F}(V)$ then the first addendum of \eqref{eq.fitzp.sum+} is convex because 
of the linearity of $L$. As the second term is quadratic and positive, $g_1$ is also convex.
 
\smallskip
(iii) As the graph of $\alpha$ is weakly closed, $g_2$ is weakly lower semi-continuous. 
As
\[
I_\alpha(v,v^* -Lv) + \langle v^*,v \rangle =I_\alpha(v,v^* -Lv) + \langle Lv,v \rangle
\qquad\forall v,v^*)\in V \!\times\! V',
\]
$g_2$ fulfills  \eqref{eq.fitzp.nonconv2} and \eqref{eq.fitzp.nonconv3}. 
\hfill$\Box$ 

\bigskip 
Next we need a lemma.

\begin{lemma} \label{lem1} 
Let $X$ be a Hausdorff topological space, $Y$ be either a reflexive space or the dual of a separable Banach space equipped with the weak star topology, and
$\varphi: X \!\times\! Y\to \erre\cup \{+\infty\}$ be lower semi-continuous. If
\begin{equation}\label{eq.fitzp.coee}
\inf_{x\in X} \; \varphi(x,y) \to +\infty
\quad\hbox{ as }\|y\|_Y \to +\infty, 
\end{equation} 
then the function $X\to \erre\cup \{+\infty\}: x\mapsto \inf_{y\in Y} \varphi(x,y)$ 
is lower semi-continuous.  
\end{lemma}  

\medskip
\noindent{\bf Proof.\/}  
This follows from classical compactness theorems, since by the next lemma
in the minimization one can confine $y$ to a weakly compact subset of $Y$. 
\hfill$\Box$

\begin{lemma} [\cite{AuEk} p.\ 120] \label{lem2} 
Let $X$ be a Hausdorff topological space, $K$ be a compact space, and
$\varphi: X \!\times\! K\to \erre\cup \{+\infty\}$ be lower semi-continuous.
The function $X\to \erre: x\mapsto \inf_{y\in K} \varphi(x,y)$ 
is then lower semi-continuous. 
\end{lemma}  

The next result extends Theorem~\ref{sumcor} to the sum of two nonlinear operators,
and generalizes Proposition~2.6 of \cite{Vi14}.
This may also be extended to finite sums of representable operators.

\begin{theorem} \label{sum}
Let $V$ be a real reflexive Banach space, $\alpha_i :V \to {\cal P}(V')$ for $i=1,2$, 
and assume that
\begin{eqnarray} 
&\exists \widetilde v\in V: \alpha_1(\widetilde v)\not=\emptyset \hbox{ and } 
\alpha_2(\widetilde v) \not=\emptyset,
\\
&g_i\in {\cal E}_w(V) \hbox{ and $g_i$ weakly represents } \alpha_i \; (i=1,2),
\\  
&\inf_{(v,v^*)\in V \!\times\! V'} \big\{g_1(v,v^*-z^*) + g_2(v,z^*)\big\} \to +\infty
\quad\hbox{ as }\|z^*\|_{V'} \to +\infty.
\label{eq.fitzp.coo} 
\end{eqnarray}
\indent
Then: 
(i) The {\rm partial inf-convolution\/} 
\begin{equation}\label{eq.fitzp.infcon}
(g_1 \!\oplus\! g_2) (v,v^*) := 
\inf_{z^*\in V'} \big\{g_1(v,v^*-z^*) + g_2(v,z^*)\big\}
\qquad\forall (v,v^*)\in V \!\times\! V'
\end{equation}
is an element of ${\cal E}_w(V)$, and represents the operator $\alpha_1 + \alpha_2$.
 
\indent 
(ii) If also $g_1,g_2\in {\cal F}(V)$, then $g_1 \!\oplus\! g_2\in {\cal F}(V)$.  
\end{theorem}

\medskip
\noindent{\bf Proof.\/}
(i) Let us  
\begin{equation}
\varphi(v,v^*,z^*) = g_1(v,v^*- z^*) + g_2(v,z^*)
\qquad\forall (v,v^*,z^*)\in V \!\times\! V'\!\times\! V'.
\end{equation} 

As this mapping is weakly lower semi-continuous, by 
\eqref{eq.fitzp.coo} and Lemma~\ref{lem1} we conclude that 
$g_1 \!\oplus\! g_2$ is weakly lower semi-continuous.
As $g_1,g_2\ge \pi$ (see \eqref{eq.pi}), it is straightforward to check that 
$g_1 \!\oplus\! g_2\ge \pi$,
and that $g_1 \!\oplus\! g_2$ represents the operator $\alpha_1 + \alpha_2$.

\smallskip
(ii) Let us define the auxiliary function
\[
\varphi(u,v^*,z^*) = g_1(u,v^*-z^*) + g_2(u,z^*)
\qquad\forall (u,v^*,z^*)\in V \!\times\! V'\!\times\! V'.
\]
If $g_1$ and $g_2$ are convex, then so is also $\varphi$. For any 
$(v_i,v^*_i,z^*_i)\in V \!\times\! V'\!\times\! V'$ ($i=1,2$) and any $\lambda\in [0,1]$, thus
\[
\varphi(\lambda(v_1,v^*_1,z^*_1) + (1-\lambda)(v_2,v^*_2,z^*_2))
\le\lambda\varphi(v_1,v^*_1,z^*_1) + (1-\lambda) \varphi(v_2,v^*_2,z^*_2). 
\]
By taking the infimum with respect to $z^*_1$ and $z^*_2$, we then get 
\[
g_1 \!\oplus\! g_2(\lambda(v_1,v^*_1)+ (1-\lambda)(v_2,v^*_2)) 
\le \lambda g_1 \!\oplus\! g_2(v_1,v^*_1)+ (1-\lambda)g_1 \!\oplus\! g_2(v_2,v^*_2);
\]
so $g_1 \!\oplus\! g_2$ is indeed convex. This yields part (ii).
\hfill$\Box$

\section{Semi-monotone operators and their representation} 
\label{sec.BEN}

\noindent 
In this section we introduce a class of nonmonotone operators
and study their variational representation, 
in the sense that we defined in Section~\ref{sec.repres}.  

\bigskip  
\noindent{\bf Semi-monotone operators.}
We shall assume throughout that $V$ is a reflexive real Banach space,
that $H$ is a real Hilbert space, and that
\begin{equation}\label{eq.BEN.spaces1}
V\subset H=H'\subset V' \text{ \ with compact and dense injections. }
\end{equation} 

If $X$ is a Banach space, by $X_s$ ($X_w$, resp.) we shall denote the same space equipped 
with the strong (weak, resp.) topology.  
Let us assume that  
\begin{eqnarray} 
&\beta(z,\cdot):V\to {\cal P}(V')\text{ is maximal monotone, }\forall z\in H, 
\label{eq.SM1}
\\ [2mm]
&\begin{split} 
&\forall (z,u)\in H\times V, \forall u^*\in \beta(z,u), 
\forall \text{ sequence $\{z_n\}$ such that $z_n\to z$ in }H,
\\
&\text{a sequence $\{u_n^*\}$ exists in $V'$ such that 
$u_n^*\in \beta(z_n,u) \;\forall n$, and $u_n^*\to u^*$ in }V'.
\label{eq.SM2}
\end{split} 
\end{eqnarray} 
This condition expresses the lower semi-continuity of the multi-valued mapping 
$\beta(\cdot,u)$, for any $u\in V$.  
This entails that $\beta(\cdot,u)$ is Borel-measurable, see e.g.\ \cite{AuFr} p.~311. 

Let us next set 
\begin{equation}\label{eq.semim}
\alpha: V\to {\cal P}(V'),
\qquad
\alpha(v) := \beta(v,v) \qquad\forall v\in V.
\end{equation}  
If \eqref{eq.SM1} and \eqref{eq.SM2} are fulfilled, we shall say that the operator 
$\alpha$ is {\it semi-monotone.\/} A differential example is displayed ahead in this section. 

The above properties entail that $\beta$ and $\alpha$ are Borel-measurable respectively on $H\times V$ and $V$.
Next we show that the graph of $\beta(\cdot,u)$ is closed for any $u\in V$.
 
\begin{lemma}\label{lem.Ken} 
If $\beta: H\times V\to {\cal P}(V')$ fulfills \eqref{eq.SM1} and \eqref{eq.SM2}, then
\begin{equation}\label{eq.SM3}
\begin{split} 
&\forall (z,u)\in H\times V, \forall \text{ sequence $\{z_n\}$ such that $z_n\to z$ in }H,
\\
&\text{if $u_n^*\in \beta(z_n,u)\;\forall n$, and $u_n^*\to u^*$ in $V'$, then }
u^*\in \beta(z,u).
\end{split} 
\end{equation}
\end{lemma}

\noindent{\bf Proof.\/}
Let two sequences $\{z_n\}$ and $\{u_n^*\}$ be as in \eqref{eq.SM3}, $w\in V$ and 
$w^*\in \beta(z,w)$. By \eqref{eq.SM2} there exists a sequence $\{w_n^*\}$ 
such that $w_n^*\in \beta(z_n,w)$ for any $n$ and $w_n^*\to w^*$ in $V'$. 
By \eqref{eq.SM1}, $\langle u_n^* - w_n^*, u-w\rangle \ge 0$ for any $n$, 
whence passing to the limit $\langle u^* - w^*, u-w\rangle \ge 0$. 
As $w\in V$ and $w^*\in \beta(z,w)$ are arbitrary, by the maximal monotonicity of 
$\beta(z,\cdot)$ we infer that $u^*\in \beta(z,u)$.
\hfill$\Box$ 

\bigskip 
\noindent{\bf Remarks.}
(i) Our definition of semi-monotonicity amends that of Kenmochi \cite{Ke},
and seems convenient for the analysis of flows; see the next sections.
At variance with \cite{Ke},
here we allow the domain of the operator $\beta$ to be a proper subset of $H\times V$.
In other terms $\beta(v,u)$ may be the empty set, for some pair $(v,u)$.
However, by \eqref{eq.SM2}, for any $v\in H$ a $u\in V$ exists such that 
$\beta(v,u)\not= \emptyset$. 
Lemma~\ref{lem.Ken} extends the analogous result at p.~220 of \cite{Ke}.

(ii) An operator $\alpha$ is called {\it generalized pseudo-monotone\/} in the sense of 
Browder-Hess \cite{BrHe} iff
\begin{equation}\label{eq.penpse}
\begin{split} 
&\forall \text{ sequence $\{(u_n, u_n^*)\}$ in }\graph(\alpha),
\\
&u_n\wto u, \;\; u_n^*\wto u^*, \;\;
\limsup_{n\to\infty} \; \langle u_n^*, u_n\rangle \le \langle u^*,u\rangle
\\
&\Rightarrow\;\; u^*\in \alpha(u), \;\; \langle u_n^*, u_n\rangle\to \langle u^*,u\rangle.
\end{split}
\end{equation}
This extends the original definition of pseudo-monotonicity of Brezis \cite{Br0}.
By mimicking the argument of \cite{Ke} p.~220, one may see that 
any semi-monotone operator $\alpha$ fulfills \eqref{eq.penpse}. 
\hfill$\Box$ 

\bigskip  
\noindent{\bf A differential example of semi-monotone operator.} 

\begin{proposition}\label{examp} 
Let $\Omega$ be a bounded domain of $\erre^N$ ($N\ge1$), $p\in{} [2,+\infty[$, and 
assume that 
\begin{eqnarray} 
&\vec\gamma(z,\cdot): \erre^N\to {\cal P}(\erre^N)
\text{ is maximal monotone, }
\label{eq.exse2}
\\
&\vec\gamma(\cdot,\vec v): \erre\to {\cal P}(\erre^N)
\text{ is lower semi-continuous (cf.\ \eqref{eq.SM2}), }
\forall \vec v\in\erre^N,
\label{eq.exse3}
\\
&\exists c_1,c_2\in\erre^+: \forall z\in\erre,\forall \vec v\in\erre^N,
\forall \vec w\in \vec \gamma(z,\vec v), \qquad
|\vec w| \le c_1 |\vec v|^{p-1} +c_2. 
\label{eq.exse4}
\end{eqnarray}
\indent
Then the operator $W^{1,p}_0(\Omega) \to {\cal P}(W^{-1,p}(\Omega)):
u\mapsto a(u) - \nabla \!\cdot\! \vec\gamma(u,\nabla u)$ is semi-monotone. 
\end{proposition}

\noindent{\bf Proof.\/}
Let us set set $V:= W^{1,p}_0(\Omega)$ and $H:= L^2(\Omega)$. 
It suffices to show the analogous of 
the property \eqref{eq.SM2} for the mapping $\vec\gamma$, that is:
\begin{equation}\label{eq.exsemi}
\begin{split} 
&\forall (z,u)\in L^2(\Omega)\times W^{1,p}_0(\Omega), 
\forall w\in \vec \gamma(z,\nabla u),  
\\
&\forall \text{ sequence $\{z_n\}$ such that $z_n\to z$ in }L^2(\Omega),
\\
&\exists \text{ a sequence $\{w_n\}$ in $L^{p'}(\Omega)$ such that } 
\\
&w_n\in\vec \gamma(z_n,\nabla u)\text{  a.e.\ in $\Omega\; \forall n$,
and $w_n\to w$ in }L^{p'}(\Omega).
\end{split} 
\end{equation} 

Notice that $z_n\to z$ in measure in $\Omega$. 
By \eqref{eq.exse3}, a sequence of measurable functions $\{w_n\}$ exists 
such that $w_n\in \vec\gamma(z_n,\nabla u)$ a.e.\ in $\Omega$ for any $n$ and 
$w_n\to w$ in measure in $\Omega$.
Moreover, $w_n\to w$ in $L^{p'}(\Omega)$, since by \eqref{eq.exse4} the sequence 
$\{|w_n|^{p'}\}$ is equi-integrable in $\Omega$. \eqref{eq.exsemi} is thus established.
\hfill$\Box$
 
\bigskip  
\noindent{\bf Representation of semi-monotone operators.} 

\begin{lemma}\label{lem.lsc} 
For any $u,\widetilde v\in V$, the mapping 
\begin{equation}\label{eq.defga} 
\gamma_{u,\widetilde v}: V\to \erre: 
v\mapsto \sup_{u^*\in \beta(v,u)} \langle u^*, v -\widetilde v\rangle
\end{equation} 
is weakly lower semi-continuous.
\end{lemma}

\noindent{\bf Proof.\/}   
Let $\{v_n\}$ be a sequence such that $v_n\wto v$ in $V$.
Let us fix any $\eps>0$, and select $u_\eps^*\in \beta(v,u)$ such that
\begin{equation}\label{eq.su}
\begin{split} 
&\langle u_\eps^*, v -\widetilde v\rangle \ge \gamma_{u,\widetilde v}(v)-\eps
\qquad\text{ if $\gamma_{u,\widetilde v}(v)< +\infty$,}
\\[2mm]
&\langle u_\eps^*, v -\widetilde v\rangle \ge 1/\eps
\qquad\qquad\quad\;\text{ otherwise.}
\end{split}
\end{equation}
We shall assume that $\gamma_{u,\widetilde v}(v)$ is finite; however this argument 
is easily extended to the case of $\gamma_{u,\widetilde v}(v) =+\infty$.

By \eqref{eq.SM2}, there exists a sequence $\{u_{\eps,n}^*\}$ in $V'$ such that 
\begin{equation}\label{eq.su=}
u_{\eps,n}^*\in \beta(v_n,u) \qquad\forall n, \qquad
u_{\eps,n}^*\to u_\eps^* \qquad\text{ in }V'. 
\end{equation}
Therefore
$\langle u_{\eps,n}^*, v_n -\widetilde v\rangle 
\to \langle u_\eps^*, v -\widetilde v\rangle$;
hence, by \eqref{eq.su}$_1$, 
\[
\liminf_{n\to \infty} \; \sup_{u^*\in \beta(v_n,u)} 
\langle u^*, v_n -\widetilde v\rangle \ge
\liminf_{n\to \infty} \; \langle u_{\eps,n}^*, v_n -\widetilde v\rangle -\eps   
\ge \langle u_\eps^*, v -\widetilde v\rangle -\eps
\ge \gamma_{u,\widetilde v}(v) - 2\eps.
\]
As $\eps$ is arbitrary, we conclude that 
\[ 
\liminf_{n\to \infty} \; \sup_{u^*\in \beta(v_n,u)} 
\langle u^*, v_n -\widetilde v\rangle 
\ge \sup_{u^*\in \beta(v,u)} \langle u^*, v -\widetilde v\rangle.
\] 

Thus $\gamma_{u,\widetilde v}$ is weakly lower semi-continuous. 
The final assertion of the lemma follows from \eqref{eq.BEN.spaces1}.
\hfill$\Box$

\bigskip
We are now able to {\it represent\/} semi-monotone operators, 
in the sense of Section~\ref{sec.repres}.

\begin{theorem}\label{BEN.param}
Let \eqref{eq.SM1}---\eqref{eq.semim} hold.
For any $z\in H$, let $f_z$ be the Fitzpatrick function of a 
maximal monotone operator $\beta(z,\cdot)$, i.e., 
\begin{equation}\label{eq.BEN.=}
f_z(v,v^*) = \sup \big\{ \langle v^*, \widetilde v\rangle - 
\langle \widetilde v^*, \widetilde v- v\rangle \!:
\widetilde v\in V, \widetilde v^* \!\in \beta(z,\widetilde v) \big\} 
\quad\;\forall (v,v^*)\in V \!\times\! V'.
\end{equation} 
Let us set
\begin{equation}\label{eq.BEN.ppp}
\varphi(v,v^*) := f_v(v,v^*)
\qquad\forall (v,v^*) \in V\times V'.
\end{equation} 
\indent
Then $\varphi\in {\cal E}_w(V)$, and it represents the operator $\alpha$.
\end{theorem}  

\noindent{\bf Proof.\/}
By the Fitzpatrick Theorem~\ref{teo.Fi}, for any $z\in H$,
\begin{equation}
\begin{split} 
&f_z(v,v^*) \ge \langle v^*,v\rangle
\qquad\forall (v,v^*)\in V \!\times\! V',
\\
&f_z(v,v^*) = \langle v^*,v\rangle
\quad\Leftrightarrow\quad v^*\in \beta(z,v).
\end{split}
\end{equation}
Setting $z=v$, we get
\begin{equation}\label{eq.BEN.a}
\begin{split} 
&\varphi(v,v^*) \ge \langle v^*,v\rangle
\qquad\forall (v,v^*)\in V \!\times\! V',
\\
&\varphi(v,v^*) = \langle v^*,v\rangle
\quad\Leftrightarrow\quad v^*\in \alpha(v).
\end{split}
\end{equation}
The function $\varphi$ thus algebraically represents the operator $\alpha$. 
By \eqref{eq.BEN.=} and \eqref{eq.BEN.ppp},
\begin{equation}\label{eq.BEN.a=}
\varphi(v,v^*) = \sup_{\widetilde v\in V} \sup_{\widetilde v^*\in \beta(v,\widetilde v)} 
\{\langle v^*, \widetilde v\rangle - \langle \widetilde v^*, \widetilde v - v\rangle \} 
\qquad\forall (v,v^*)\in V \!\times\! V'.
\end{equation} 

By Lemma~\ref{lem.lsc}, for any $\widetilde v\in V$ the mapping 
\[
(v,v^*)\mapsto \sup_{\widetilde v^*\in \beta(v,\widetilde v)}
\{\langle v^*, \widetilde v\rangle - \langle \widetilde v^*, \widetilde v - v\rangle \}
\]
is weakly lower semi-continuous on $V \!\times\! V'$. 
The same then holds for its supremum as $\widetilde v$ ranges in $V$, that is, for $\varphi$; 
thus $\varphi\in {\cal E}_w(V)$.
By \eqref{eq.BEN.a}, we conclude that $\varphi$ represents $\alpha$. 
\hfill$\Box$

\bigskip  
Next we extend Theorem~\ref{BEN.param} to general representative functions.
Let us first denote by $\lsc\varphi$ the lower semi-continuous 
{\it regularized function\/} of any $\varphi:V\!\times\! V'\to \erre\cup \{+\infty\}$
w.r.t.\ the weak topology, 
namely the pointwise supremum of the weakly continuous minorants of $\varphi$.

\begin{proposition}\label{BEN.param=} 
Let \eqref{eq.SM1}---\eqref{eq.semim} hold, assume that
\begin{eqnarray} 
&g_z\text{ weakly represents the maximal monotone operator }
\beta(z,\cdot), \forall z\in H,
\label{eq.BEN.00} 
\end{eqnarray}  
and set
\begin{equation}\label{eq.BEN.pp+}
\varphi(v,v^*) := g_v(v,v^*)
\qquad\forall (v,v^*) \in V\times V'.
\end{equation}  
\indent
Then $\lsc\varphi\in {\cal E}_w(V)$, and it represents the operator $\alpha$.  
\end{proposition}

\noindent{\bf Proof.\/} 
For any $z\in H$, let us denote by $f_z$ the Fitzpatrick function of $\beta(z,\cdot)$, 
and set $\psi(v,v^*) := f_v(v,v^*)$ for any $(v,v^*) \in V\times V'$.
By Theorem~\ref{BEN.param} $\psi\in {\cal E}_w(V)$, and this function 
represents the operator $\alpha$.
By Lemma~\ref{lem.dis=}, $f_z\le g_z$ for any $z\in H$; hence 
$\psi \le \lsc\varphi \le\varphi$. 
As both $\psi$ and $\varphi$ represent the operator $\alpha$, 
the same holds for $\lsc\varphi \in{\cal E}_w(V)$. 
\hfill$\Box$  

\medskip
\noindent{\bf An example.} 

\begin{proposition}\label{eq.fitp.ellpse} 
Let $\Omega$ be a bounded domain of $\erre^N$ ($N\ge 1$) of Lipschitz class, 
$p\in{} [2,+\infty[$, set $V:= W^{1,p}_0(\Omega)$, $H:= L^2(\Omega)$, 
define $\Lambda$ as in \eqref{eq.fitp.du}, and let $\vec\gamma_z$ 
($= \vec\gamma(z,\cdot)$) fulfill \eqref{eq.exse2}--\eqref{eq.exsemi}.
Assume that, for any $z\in \erre$,
\begin{eqnarray}  
&\erre\to \erre\cup \{+\infty\}: z\mapsto g_z(v,v')
\text{ is Borel measurable, }\forall (v,v')\in \erre^N\!\times\! \erre^N,
\label{eq.BEN.00=}
\\
&g_z\text{ represents the maximal monotone operator }
\vec\gamma_z: \erre^N\to {\cal P}(\erre^N), 
\label{eq.BEN.00-}  
\end{eqnarray}
and set
\begin{equation}\label{eq.fitzp.a2=}
\varphi(v,v^*) = \int_\Omega g_{v(x)}(\nabla v(x), \nabla \Lambda v^*(x)) \, dx
\qquad\forall (v,v^*)\in V \!\times\! V',
\end{equation}  
\indent
Then $\lsc\varphi\in {\cal E}_w(V)$, and it represents the semi-monotone operator 
\begin{equation}\label{eq.fitzp.a2x}
V\to {\cal P}(V'): v\mapsto - \nabla \cdot \vec\gamma(v,\nabla v). 
\end{equation} 
\indent
If $g_z$ is the Fitzpatrick function of $\vec\gamma_z$ for any $z\in \erre$, 
then $\varphi\in {\cal E}_w(V)$.
\end{proposition} 

\noindent{\bf Proof.\/} 
By Proposition~\ref{eq.fitp.ell}, for any $z\in H$ the mapping
\begin{equation}\label{eq.fitzp.a2+}
\psi_z(v,v^*) = \int_\Omega g_{z(x)}(\nabla v(x), \nabla \Lambda v^*(x)) \, dx
\qquad\forall (v,v^*)\in V \!\times\! V'
\end{equation}  
represents the maximal monotone operator 
$\beta(z,\cdot): V\to {\cal P}(V'): v\mapsto - \nabla \cdot \vec\gamma(z,\nabla v)$.
The thesis then follows from Theorem~\ref{BEN.param} and Proposition~\ref{BEN.param=}.
\hfill$\Box$

\section{Time-dependence} 
\label{sec.time} 

In this section we extend the previous results to spaces of time-dependent functions,
as a step towards representing semi-monotone flows in the next section.

Let us assume \eqref{eq.BEN.spaces1}, fix any $p\in [2,+\infty[$, and notice that 
\begin{equation}\label{eq.BEN.spaces1-} 
\begin{split}
&L^p(0,T;V) \subset L^2(0,T;H) = L^2(0,T;H)' 
\subset L^p(0,T;V)' \simeq L^{p'}((0,T;V')
\\
&\text{with (noncompact) continuous and dense injections. } 
\end{split} 
\end{equation}
Let us also define
\begin{equation}\label{eq.BEN.spaces2}
{\cal V} := \big\{v\in L^p(0,T;V): D_tv \in L^{p'}(0,T;V')\big\}.
\end{equation}
Equipped with the graph norm, this is a Banach space and can be identified with 
a subspace of $C^0([0,T];H)$. 
The space ${\cal V}$ arises in a natural way in the analysis of parabolic equations,
since, because of \eqref{eq.BEN.spaces1} and
by the classical Aubin's lemma \cite{Au},
\begin{equation}\label{eq.BEN.comp}
{\cal V} \subset L^2(0,T;H) \text{ \ with compact and dense injection. } 
\end{equation} 

Let \eqref{eq.SM1}---\eqref{eq.semim} also hold, and assume that
\begin{equation}\label{eq.BEN.ineq} 
\begin{split}  
&\exists A,B>0: \forall z\in H,\forall v\in V, \forall \xi\in \beta(z,v)\qquad
\\
&\|\xi\|_{V'} \le A \max\{\|v\|_V^{p-1}, \|z\|_H^{p-1}\} +B,
\end{split} 
\end{equation} 
whence, by the boundedness of the injection $V\to H$,
for a suitable constant $\widetilde A>0$,
\[  
\|\xi\|_{V'} \le \widetilde A\|v\|_V^{p-1} +B \qquad\forall \xi\in\alpha(v), \forall v\in V. 
\]

Therefore $\beta$ and $\alpha$ canonically determine two global-in-time operators:
\begin{eqnarray}  
&\begin{split}\label{eq.BEN.sum1}
&\widehat\beta: L^p(0,T;H \!\times\! V)\to {\cal P}(L^{p'}(0,T;V')):
\\
&[\widehat\beta(z,v)](t) = \beta(z(t),v(t)) 
\\
&\text{for a.e.\ }t\in{} ]0,T[, \forall (z,v)\in L^p(0,T;H \!\times\! V), 
\end{split}
\\
&\begin{split}\label{eq.BEN.sum2}
&\widehat\alpha: L^p(0,T;V)\to {\cal P}(L^{p'}(0,T;V')):
\\
&\widehat\alpha(v) = \widehat\beta(v,v)
\quad\text{ a.e.\ in }]0,T[,\forall v\in L^p(0,T;V),
\end{split}
\\
&\begin{split}\label{eq.BEN.sum3}
&\|\xi\|_{L^{p'}(0,T;V')} \le \widetilde A\|v\|_{L^p(0,T;V)}^{p-1} + BT^{1/p'} 
\\
&\forall \xi\in \widehat\alpha(v), \forall v\in VL^p(0,T;V). 
\end{split}
\end{eqnarray} 

Let us now point out a simple property, which plays a role in representing
operators acting on time-dependent functions. 

\begin{lemma}\label{fundlem}
Whenever $\varphi$ represents an operator $\alpha: V\to {\cal P}(V')$ 
(in the sense of Section~\ref{sec.repres}), the corresponding integral functional 
$\Phi: (v,v^*) \mapsto  \int_0^T \varphi(v(t),v^*(t)) \, dt$ algebraically
represents the corresponding operator 
$\widehat\alpha: L^p(0,T;V) \to {\cal P}(L^{p'}(0,T;V'))$.
\end{lemma}

\noindent{\bf Proof.\/}
As $\varphi\ge \pi$, 
\begin{equation}\label{eq.BEN.int} 
\int_0^T [\varphi(v,v^*) - \langle v^*, v\rangle] \, dt \ge 0
\qquad\forall (v,v^*)\in L^p(0,T;V) \!\times\! L^{p'}(0,T;V').
\end{equation}
Moreover, $\int_0^T \varphi(v,v^*) \, dt= \int_0^T \langle v^*, v\rangle \, dt$ for some 
$(v,v^*)\in L^p(0,T;V) \!\times\! L^{p'}(0,T;V')$
only if the nonnegative integrand of \eqref{eq.BEN.int} vanishes.
As $\varphi\ge \pi$, this means that $\varphi(v,v^*) = \pi(v,v^*)$ a.e.\ in $]0,T[$. 
As $\varphi$ represents $\alpha$, we conclude that $v^*\in \alpha(v)$ a.e.\ in $]0,T[$,
that is, $v^*\in \widehat\alpha(v)$.
\hfill$\Box$
 
\begin{proposition}\label{intrep=}
Let $p\in [2,+\infty[$, let \eqref{eq.SM1}---\eqref{eq.semim},
\eqref{eq.BEN.ineq}--\eqref{eq.BEN.sum2} hold, assume that
\begin{eqnarray}  
&\begin{split}
&H\to \erre\cup \{+\infty\}: z\mapsto g_z(v,v^*)
\text{ is Borel measurable, }\forall (v,v^*)\in V\!\times\! V', 
\end{split}
\label{eq.time.1}
\\ 
&g_z\in {\cal F}(V), \text{ and it represents }
\beta(z,\cdot): V\to {\cal P}(V'), \forall z\in H,
\label{eq.time.2}
\end{eqnarray}   
and set
\begin{equation}\label{eq.BEN.t}
\begin{split}
&G_z(v,v^*) :=  \int_0^T g_{z(t)}(v(t),v^*(t)) \, dt 
\\
&\forall (z,v)\in L^p(0,T;H \!\times\! V), \forall v^*\in L^{p'}(0,T;V').
\end{split} 
\end{equation}  
\indent
Then:
(i) For any $z\in L^p(0,T;H)$, $G_z\in {\cal F}(L^p(0,T;V))$ and it represents 
$\widehat\beta(z,\cdot)$.

\indent
(ii) If $g_z$ is the Fitzpatrick function of $\beta(z,\cdot)$ for any $z\in H$ 
(see \eqref{eq.fitzp.Fitzfun}), 
then $G_z$ is the Fitzpatrick function of $\widehat\beta(z,\cdot)$ for any $z\in L^p(0,T;H)$.
\end{proposition}

\noindent{\bf Proof.\/}
(i) Let us fix any $z\in L^p(0,T;H)$, and define the mapping 
\[ 
\theta_z(v,v^*) := G_z(v,v^*) - \!\! \int_0^T \!\langle v^*, v\rangle \, dt  
\qquad\forall v\in L^p(0,T;V), \forall v^*\in L^{p'}(0,T;V').
\]
As $g_z$ is a representative function, $\theta_z\ge 0$.
By Fatou's lemma, by the semi-continuity of $g_z$ 
and by the continuity of the duality pairing, $\theta_z$ is lower semi-continuous on 
$V \!\times\! V'$.
The same then holds for $G_z$ on $L^p(0,T;V) \!\times\! L^{p'}(0,T;V')$. 
By Lemma~\ref{fundlem}, we conclude that $G_z$ represents $\widehat\beta(z,\cdot)$.

\smallskip
(ii) For any $z\in L^p(0,T;H)$ the operator $\widehat\beta(z,\cdot)$ 
is maximal monotone. By \eqref{eq.fitzp.Fitzfun}, its Fitzpatrick function reads
\begin{equation}\label{eq.BEN.+}
\begin{split}
&F_z(v,v^*) = \sup \bigg\{ \int_0^T [\langle v^*, \widetilde v\rangle - 
\langle \widetilde v^*, \widetilde v- v\rangle] \, dt:
\\
&(\widetilde v,\widetilde v^*)\in L^p(0,T;V) \!\times\! L^{p'}(0,T;V'), 
\widetilde v^* \!\in \widehat\beta(z,\widetilde v) \bigg\},
\\
&\forall v\in L^p(0,T;V), \forall v^*\in L^{p'}(0,T;V').
\end{split}
\end{equation}

If $v,v^*,\widetilde v,\widetilde v^*$ are all assumed to be simple functions of time
(i.e., measurable and with finite range), then
the supremum of the integrals coincides with the integral of the supremum,
namely with the integral of the Fitzpatrick function $g_{z(\cdot)}$ of 
$\beta(z(\cdot),\cdot)$; that is, $F_z = G_z$.
As simple functions are dense in the space of integrable functions, 
the same equality then holds for the latter class, too.    
\hfill$\Box$ 

\bigskip
Next we show that $\widehat\alpha$ is semi-monotone, and represent it.

\begin{theorem}\label{timeint}
Let $p\in [2,+\infty[$, let \eqref{eq.SM1}---\eqref{eq.semim},
\eqref{eq.BEN.ineq}--\eqref{eq.BEN.sum2} hold, and 
let $f_z$ be the Fitzpatrick function of $\beta(z,\cdot)$ for any $z\in H$
(see \eqref{eq.fitzp.Fitzfun}). Then:

\indent
(i) The operator $\widehat\alpha: L^p(0,T;V)\to {\cal P}(L^{p'}(0,T;V'))$ 
is semi-monotone.

\indent 
(ii) Let us set $\varphi(v,v^*) := f_v(v,v^*)$ for any $(v,v^*) \in V\times V'$ and
\begin{equation}\label{eq.BEN.pp=}
\begin{split} 
&\Phi(v,v^*) := \left\{\! \begin{split}
&\int_0^T \varphi(v(t),v^*(t)) \, dt
\qquad\forall (v,v^*) \in {\cal V} \times (L^{p'}(0,T;V'),
\\
&+\infty
\qquad\;\forall (v,v^*) \in 
[L^p(0,T;V) \setminus {\cal V}] \times (L^{p'}(0,T;V').
\end{split}\right. 
\end{split}
\end{equation}
Then $\Phi\in {\cal E}_w(L^p(0,T;V))$, and it represents the operator $\widehat\alpha$.

\indent
(iii) Assume that $g_z$ is any function that fulfills \eqref{eq.time.1} and \eqref{eq.time.2}, 
set $\bar \varphi(v,v^*) := g_v(v,v^*)$ for any $(v,v^*) \in V\times V'$, 
and define $\widetilde\Phi$ as in \eqref{eq.BEN.pp=},
with $\widetilde\varphi = {\rm sc_s^-\;} \bar \varphi$ 
(the lower semi-continuous regularized function w.r.t.\ the strong topology) 
in place of $\varphi$.  
Then $\widetilde\Phi\in {\cal E}_s(L^p(0,T;V))$, and it represents $\widehat\alpha$.
\end{theorem} 

\noindent{\bf Proof.\/} 
(i) For any $z\in L^p(0,T;H)$, $\widehat\beta(z,\cdot)$ is maximal monotone.  

To prove the analogous property of \eqref{eq.SM2} for $\widehat\beta$, we must show that,
for any $(z,u)\in L^p(0,T;H\times V)$, 
for any $u^*\in \widehat\beta(z,u)$, and
for any sequence $\{z_n\}$ such that 
\begin{equation}\label{eq.BEN.r}
z_n\to z \qquad\text{ in }L^p(0,T;H),
\end{equation}
there exists a sequence $\{u_n^*\}$ in $L^{p'}(0,T;V')$ such that 
\begin{equation}\label{eq.BEN.s}
u_n^*\in \widehat\beta(z_n,u) \;\forall n, 
\qquad
u_n^*\to u^* \qquad\text{ in }L^{p'}(0,T;V').
\end{equation}
 
By \eqref{eq.BEN.r}, up to extracting a subsequence, $z_n\to z$ in $H$ a.e.\ in $]0,T[$.
By \eqref{eq.SM2}, for a.e.\ $t\in{} ]0,T[$ then there exists a $t$-measurable sequence 
$\{u_n^*(t)\}$ in $V'$ such that $u_n^*(t)\in \widehat\beta(z_n(t),u(t))$ for any $n$ and 
$u_n^*(t)\to u^*(t)$ in $V'$. We can also assume that for any $n$ the function 
$u_n^*:{}]0,T[{}\to V'$ is measurable. 

Notice that, by \eqref{eq.BEN.ineq}, $u^*\in L^{p'}(0,T;V')$ and the sequence $\{u_n^*\}$ 
is bounded in this space. Moreover, for any measurable subset $M$ of $]0,T[$,
\begin{equation}\label{eq.equi} 
\int_M \|u_n^*\|_{V'}^{p'} \, dt \le  
A\! \int_M \|z_n\|_{V'}^p \, dt + B |M| 
\quad\forall n, \forall \text{ measurable }M\subset{} ]0,T[.
\end{equation} 
This entails the following property of equi-integrability:
\begin{equation}\label{eq.equi=} 
\sup_{n\in \enne} \int_M \|u_n^*\|_{V'}^{p'} \, dt \to 0
\quad\text{ as }|M|\to 0.
\end{equation} 
Hence $u_n^*\to u^*$ in $L^{p'}(0,T;V')$.
The analogous of property \eqref{eq.SM2} is thus established for $\widehat\beta$. 
Therefore $\widehat\alpha$ is semi-monotone.

\smallskip
(ii) Let us define $F_z$ as in \eqref{eq.BEN.t}  for any $z\in H$, 
here with $f_z$ in place of $g_z$.
By Proposition~\ref{intrep=}, $F_z\in {\cal F}(L^p(0,T;V))$ 
and it is the Fitzpatrick function of $\widehat\beta(z,\cdot)$ for any $z\in L^p(0,T;H)$. 
By Theorem~\ref{BEN.param} (here with $L^p(0,T;V)$ in place of $V$), we 
then conclude that $\Phi\in {\cal E}_w(L^p(0,T;V))$, and it represents $\widehat\alpha$.  

\smallskip
(iii) Next let us assume that $g_z$ fulfills \eqref{eq.time.1} and \eqref{eq.time.2},
and define $\widetilde\Phi$ as in \eqref{eq.BEN.pp=} with $\widetilde\varphi = 
{\rm sc_s^-\;}\bar\varphi$ in place of $\varphi$.
By Lemma~\ref{fundlem}, as ${\rm sc_s^-\;} \bar\varphi$ algebraically represents 
$\alpha$ in $V$,
$\widetilde\Phi$ algebraically represents $\widehat\alpha$ in $L^p(0,T;V)$.

We are left with proving that $\widetilde\Phi$ is lower semicontinuous.
If $(u_n,u_n^*)\to (u,u^*)$ in $L^p(0,T;V) \!\times L^{p'}(0,T;V')$ then, up to extracting a
subsequence, $(u_n,u_n^*)\to (u,u^*)$ in $V\!\times V'$ a.e.\ in $]0,T[$.
Hence $\liminf_{n\to \infty} \widetilde\varphi(u_n,u_n^*)\ge \widetilde\varphi(u,u^*)$ a.e.\ in 
$]0,T[$. As by \eqref{eq.BEN.ineq} the sequence $\{\widetilde\varphi(u_n,u_n^*)\}$ 
is equi-integrable in $]0,T[$, we infer that 
$\liminf_{n\to \infty} \widetilde\Phi(u_n,u_n^*)\ge \widetilde\Phi(u,u^*)$.
\hfill$\Box$ 

\begin{proposition}\label{pse}
Let $p\in [2,+\infty[$ and let \eqref{eq.SM1}---\eqref{eq.semim},
\eqref{eq.BEN.ineq}--\eqref{eq.BEN.sum2} hold. Let $u^0\in H$ and set  
\begin{eqnarray}
&{\cal V}_{u^0} :=  \{v\in{\cal V}_{u^0}: v(0)= u^0 \},
\label{eq.BEN.qs}
\\
&R(z,v) :=  
\left\{\begin{split}
&D_t v+ \widehat\beta(z,v) 
\qquad\forall z\in L^p(0,T;H), \forall v\in {\cal V}_{u^0},
\\[2mm]
&\emptyset
\qquad\forall z\in L^p(0,T;H), \forall v\in L^p(0,T;V)\setminus {\cal V}_{u^0}.
\end{split}\right.
\label{eq.BEN.qq}
\end{eqnarray} 
\indent
Then $R(z,\cdot)$ is maximal monotone, and
\begin{eqnarray}\label{eq.comstab.mmm+}
&\begin{split} 
&S: L^p(0,T;V) \to {\cal P}(L^{p'}(0,T;V')): v\mapsto R(v,v) 
\text{ \ is semi-monotone.}
\end{split}
\end{eqnarray}
\end{proposition}

\noindent{\bf Proof.\/}
For any $z\in L^p(0,T;H)$, the operator 
$D_t+ \widehat\beta(z,\cdot) : L^p(0,T;V) \to {\cal P}(L^{p'}(0,T;V'))$ with domain 
${\cal V}_{u^0}$ is monotone, since so is $\widehat\beta(z,\cdot)$ and
\[
\int_0^T \langle Dv(t),v(t) \rangle \, dt = {1\over2}\|v(T)\|_H^2 \ge 0
\qquad \forall v\in {\cal V}_{u^0}- {\cal V}_{u^0}.
\]

Next let $\Sigma$ be the duality mapping of the space $V$.
By a classical theory, see e.g.\ \cite{Bar}, \cite{Br1}, for any $z\in L^p(0,T;H)$ and any 
$u^*\in L^{p'}(0,T;V')$, there exists a (unique) solution $u\in {\cal V}_{u^0}$ of the 
equation
\[
D_tu + \beta(z,u) +\Sigma(u) = u^*
\qquad \text{ in $V'$, a.e.\ in }{}]0,T[.
\]  

By the Minty-Browder theorem, we conclude that $R(z,\cdot)$ is maximal monotone. 
In the proof of Theorem~\ref{timeint} we saw that the operator $\widehat\beta$ fulfills the analogous of property \eqref{eq.SM2}. The same then holds for $R$, and 
therefore $S$ is semi-monotone. 
\hfill$\Box$

\section{Representation of semi-monotone flows} 
\label{sec.psflow}

\noindent 
In this section we extend the Brezis-Ekeland-Nayroles (BEN) principle to the flow of a semi-monotone operator $\alpha: V\to {\cal P}(V')$ defined as in \eqref{eq.semim}.  
 
Let us fix any $u^*\in L^{p'}(0,T;V')$, any $u^0\in H$, 
define ${\cal V}_{u^0}$ as in \eqref{eq.BEN.qs}, and consider the Cauchy problem  
\begin{equation}\label{CauPb1}  
u\in {\cal V}_{u^0},
\qquad
D_tu + \alpha(u) \ni u^* \qquad \text{ in $V'$, a.e.\ in }{}]0,T[. 
\end{equation} 
Defining $\widehat\alpha: L^p(0,T;V) \to {\cal P}(L^{p'}(0,T;V'))$ 
as in \eqref{eq.BEN.sum2},
this problem can equivalently be formulated globally in time as follows:
\begin{equation}\label{CauPb2}  
u\in {\cal V}_{u^0},
\qquad
D_tu + \widehat\alpha(u) \ni u^* \qquad\text{in }L^{p'}(0,T;V'). 
\end{equation} 

Next we deal with the representation of the operator 
$D_t+ \widehat\alpha: L^p(0,T;V) \to {\cal P}(L^{p'}(0,T;V'))$
with domain ${\cal V}_{u^0}$
(this coincides with the operator of \eqref{eq.comstab.mmm+}).

\begin{theorem}(Extended BEN principle) \label{prop.repr}
Let $p\in [2,+\infty[$, let us assume that 
\eqref{eq.SM1}---\eqref{eq.semim}, \eqref{eq.BEN.ineq} hold.  
Let $f_z$ be the Fitzpatrick function of $\beta(z,\cdot)$ for any $z\in H$
(see \eqref{eq.fitzp.Fitzfun}), and set 
\begin{gather}
\varphi(v,v^*) := f_v(v,v^*) \qquad \forall (v,v^*) \in V\times V',
\\
\begin{split} \label{eq.comstab.BENfun+} 
&\Gamma(v,v^*) := \int_0^T [\varphi(v,v^* -D_tv) - \langle v^*- D_tv,v\rangle] \, dt   
\\
&\qquad\quad\;\; = \int_0^T [\varphi(v,v^* -D_tv) - \langle v^*,v\rangle] \, dt 
+ {1\over2} \|v(T)\|_H^2 - {1\over2} \|u^0\|_H^2  
\\
& \qquad\qquad \qquad \qquad \qquad\qquad\qquad\qquad
\forall (v,v^*)\in {\cal V}_{u^0} \!\times\! L^{p'}(0,T;V'), 
\\
&\Gamma(v,v^*) := +\infty \quad\;\;
\text{for any other }(v,v^*)\in L^p(0,T;V) \times L^{p'}(0,T;V').
\end{split}
\end{gather} 
\indent
Then:  
(i) $\Gamma\in {\cal E}_w(L^p(0,T;V))$, 
and it represents the operator $D_t + \widehat\alpha$. 
 
\indent 
(ii) Let instead $f_z$ be any representative function of $\beta(z,\cdot)$ for any $z\in H$
that fulfills \eqref{eq.time.1}, and define $\widetilde\Phi$ as in 
\eqref{eq.BEN.pp=} with $\widetilde\varphi = {\rm sc_s^-\;} \varphi$ in place of $\varphi$.  
Then $\widetilde\Phi\in {\cal E}_s(L^p(0,T;V))$, 
and it represents the operator $\widehat\alpha$.

\indent
(iii) In either case, the Cauchy problem \eqref{CauPb2} is equivalent to the following 
problem of {\it null-minimization:\/}
\footnote{ We refer to problems of this sort, in which the minimum is prescribed to vanish,
as {\it null-minimization problems.\/}
}
\begin{equation}\label{eq.BEN.min+} 
u\in {\cal V}_{u^0},
\qquad
\Gamma(u,u^*) = 0 \;\; \Big(\! =\inf_{{\cal V}_{u^0}} \Gamma(\cdot,u^*) \Big). 
\end{equation} 
\end{theorem}
 
\noindent{\bf Proof.\/}  
(i) Let us first define the functional $\Phi$ as in \eqref{eq.BEN.pp=};
by part (ii) of Theorem~\ref{timeint}, $\Phi$ is weakly lower semicontinuous and represents 
$\widehat\alpha$. Notice that 
\[
\Gamma(v,v^*) = \Phi(v,v^* - D_tv ) - \int_0^T \langle v^*- D_tv,v\rangle \, dt
\qquad\forall (v,v^*)\in {\cal V}_{u^0} \!\times\! L^{p'}(0,T;V').
\]

Part (i) then stems from Theorem~\ref{sum}, here applied with 
$\alpha_1 = \widehat\alpha$, $\alpha_2 = D_t$ and with the space $L^p(0,T;V)$ 
in place of $V$.

Part (ii) is proved via a similar argument, defining $\widetilde\Phi$ as in 
\eqref{eq.BEN.pp=} with $\widetilde\varphi = {\rm sc_s^-\;} \varphi$ in place of $\varphi$, 
and referring to part (iii) of Theorem~\ref{timeint}.

Part (iii) follows from parts (i) and (ii).
\hfill$\Box$ 

\bigskip
\noindent{\bf Modified set-up.}
Next we amend the functional set-up \eqref{eq.BEN.spaces1-}, \eqref{eq.BEN.spaces2},
in order to improve the robustness of the BEN-type minimization principle (see ahead). 
Let us first define the measure
\begin{equation}\label{eq.comstab.mu}
\mu(A) = \int_A (T-t) \, dt \quad\forall A\in {\cal L}(0,T),
\quad\text{ i.e. }\quad
d\mu(t) = (T-t) dt.
\end{equation}
For any reflexive Banach space $X$ and any $p\in{} ]1,+\infty[$, 
let us introduce the spaces
\begin{equation}\label{eq.comstab.spaces4} 
\begin{split} 
&L^p_\mu(0,T;X) := \big\{\mu\text{-measurable } v:{} ]0,T[{}\to X:
\|v(t)\|_X \in L^p_\mu(0,T) \big\}, 
\\ 
&{\cal V}_\mu := \big\{v\in L^p_\mu(0,T;V): D_tv \in L^{p'}_\mu(0,T;V')\big\},
\\ 
\end{split} 
\end{equation}
and equip them with the respective graph norm.  
These are reflexive Banach spaces, and for any $p\in{} [2,+\infty[$
\begin{equation} 
L^p_\mu(0,T;V) \subset L^2_\mu(0,T;H) = L^2_\mu(0,T;H)' 
\subset L^{p'}_\mu(0,T;V)' = L^p_\mu(0,T;V'),
\end{equation} 
with continuous and dense injections.  

\bigskip
\noindent{\bf Time-integrated extended BEN principle.}
Let us fix any $p\in{} [2,+\infty[$, any $u^0\in H$, and set
\begin{equation}\label{eq.comstab.spaces3}
{\cal V}_{\mu,u^0} := \big\{v\in L^p_\mu(0,T;V): 
D_tv \in L^{p'}_\mu(0,T;V'), v(0) = u^0\big\} \subset {\cal V}.
\end{equation}
($v(0)$ is meaningful since ${\cal V}_\mu \subset C^0([0,T];H)$.) 
Notice that
\begin{equation}\label{eq.comstab.double} 
\begin{split}
&\int_0^T \langle D_tv,v\rangle \, d\mu(t) 
\overset{\eqref{eq.comstab.mu}}{=} 
\int_0^T \! d\tau\! \int_0^\tau \! \langle D_tv,v\rangle \, dt 
\\
&= {1\over2}\! \int_0^T\! d\tau\! \int_0^\tau \! D_t\big(\|v(t)\|_H^2\big) \, dt 
= {1\over2} \!\int_0^T \! \|v(\tau)\|_H^2 \, d\tau 
\qquad\forall v\in {\cal V}_\mu, 
\end{split}
\end{equation}
so that 
\begin{equation}
D_t:{\cal V}_{\mu,u^0}(\subset L^p_\mu(0,T;V)) \to L^{p'}_\mu(0,T;V')
\text{ is monotone. }
\end{equation}

By double time-integration, let us define the nonnegative functional  
\begin{equation} \label{eq.comstab.BENfun'} 
\begin{split}
&\Psi(v,v^*) 
:= \!\int_0^T \!\!\! d\tau\! \int_0^\tau \!\! \big[\varphi(v,v^* -D_tv) 
- \langle v^*-D_tv,v\rangle\big] \, dt  
\\
&\qquad\quad\overset{\eqref{eq.comstab.double}}{=}
\!\!\!\! \int_0^T \!\! \big[\varphi(v,v^* -D_tv) - \langle v^*,v\rangle\big] \, d\mu(t)  
+{1\over2}\! \int_0^T \!\!\! \|v(t)\|_H^2 \, dt - {T\over2} \|u^0\|_H^2
\\
&\qquad\qquad\qquad\qquad\qquad\qquad\qquad\qquad\quad
\qquad\forall (v,v^*)\in {\cal V}_{\mu,u^0} \!\times\! L^{p'}_\mu(0,T;V'), 
\\
&\Psi(v,v^*) 
:= +\infty \qquad\qquad
\text{for any other }(v,v^*)\in L^p(0,T;V) \times L^{p'}(0,T;V').
\end{split}
\end{equation} 

\begin{theorem}(Time-integrated extended BEN principle) \label{prop.repr'} 
Let $p\in [2,+\infty[$, let \eqref{eq.SM1}---\eqref{eq.semim},
\eqref{eq.BEN.ineq}--\eqref{eq.BEN.sum2} hold.
Let $f_z$ be the Fitzpatrick function of $\beta(z,\cdot)$ for any $z\in H$
(see \eqref{eq.fitzp.Fitzfun}), 
set $\varphi(v,v^*) := f_v(v,v^*)$ for any $(v,v^*) \in V\times V'$,
and define $\Psi$ as in \eqref{eq.comstab.BENfun'}.  
\\
\indent 
Then: 
(i) $\Psi\in {\cal E}_w(L^p_\mu(0,T;V))$,  
and it represents the operator $D_t + \widehat\alpha$. 
 
\indent
(ii) Let instead $f_z$ be any representative function of $\beta(z,\cdot)$ for any $z\in H$
that fulfills \eqref{eq.time.1}, and define $\Psi$ as in 
\eqref{eq.comstab.BENfun'} with $\widetilde\varphi = {\rm sc_s^-\;} \varphi$ 
in place of $\varphi$.  
Then $\Psi\in {\cal E}_s(L^p_\mu(0,T;V))$, and it strongly represents the operator 
$\widehat\alpha$. 

\indent
(iii) In either case, the Cauchy problem \eqref{CauPb2} is equivalent to the following 
null-minimization problem:
\begin{equation}\label{eq.BEN.min++} 
u\in {\cal V}_{\mu,u^0},
\qquad
\Psi(u,u^*) = 0 \;\; 
\Big(\! =\inf_{{\cal V}_{\mu,u^0}} \Psi(\cdot,u^*) \Big). 
\end{equation}
\end{theorem} 

\noindent{\bf Proof.\/} 
This follows from Proposition~\ref{prop.repr}.
For part (iii), notice that ${\cal V}\cap {\cal V}_\mu$ in dense in ${\cal V}_\mu$,
and the further integration in $]0,T[$ maps ${\cal E}_s(L^p(0,T; V))$ to 
${\cal E}_s(L^p_\mu(0,T; V))$.
\hfill$\Box$   

\begin{remarks}\rm 
(i) By comparing the functionals \eqref{eq.comstab.BENfun+} and 
\eqref{eq.comstab.BENfun'}, the Reader will notice that in \eqref{eq.comstab.BENfun'} 
${1\over2} \|v(T)\|_H^2$ is replaced by ${1\over2} \int_0^T \!\! \|v(t)\|_H^2 \, dt$. 
The latter term is continuous w.r.t. to the weak convergence in ${\cal V}_\mu$,
at variance with the former.
This is the main reason why we introduced the further integration in time. 

(ii) The above formulation can be extended in several ways. For instance, the operator 
$\beta$, and then also $\alpha$, might explicitly depend on time. 
\hfill$\Box$  
\end{remarks}

\section{Structural compactness, structural stability and nonlinear weak topology}  
\label{sec.comstab}

\noindent
In this section we illustrate the notions of structural compactness and stability, 
and apply them to semi-monotone flows, here reformulated as null-minimization principles. 
We also introduce what we shall refer to as the {\it nonlinear weak topology\/} of 
$V \!\times\! V'$, which will be used in Theorem~\ref{eq.par.stru}. 

\bigskip
\noindent{\bf Structural compactness and stability of minimization problems.}  
First, we illustrate these notions in an abstract set-up. Let $X$ be a topological space
and ${\cal G}$ be a family of functionals $X\to \erre\cup\{+\infty\}$,
equipped with a suitable notion of convergence. 
We shall use the following terminology:

(i) the problem of minimizing these functionals will be called {\it structurally compact\/} 
if the family ${\cal G}$ is sequentially compact, and
the corresponding minimizers range in a sequentially relatively compact subset of $X$. 

(ii) the minimization problem will be said {\it structurally stable\/} if 
\begin{equation}\label{eq.comstab.sta0}
\left\{\begin{split} 
&u_n\to u \quad\hbox{ in }X
\\ 
&\Phi_n \to \Phi\quad\hbox{ in }{\cal G}
\\ 
&\Phi_n(u_n) - \inf \Phi_n \to 0
\end{split} 
\qquad\Rightarrow\quad
\Phi(u) = \inf \Phi. \right.
\end{equation}  
The definition (i) is clearly instrumental to (ii).
This also applies to null-minimiz\-ation problems.

The selection of the notion of convergence in ${\cal G}$ is crucial. 
Structural compactness and stability are in competition:
the convergence must be sufficiently weak in order to allow for sequential compactness,
and at the same time it must be so strong to provide stability
(essentially, passage to the limit in the perturbed minimization problem). 
We shall see that a notion of $\Gamma$-convergence is especially appropriate for 
this problem; here we refer to the (typically stationary) ordinary $\Gamma$-convergence,
as well as the evolutionary $\Gamma$-convergence of Section~\ref{sec.evol}.
$\Gamma$-convergence looks more convenient than other notions of variational convergence
like the Mosco-convergence of \cite{Mo},
namely the simultaneous weak and strong $\Gamma$-convergence to a same function, 
see e.g.\ \cite{At},\cite{Mo}. 
On the other hand the need of compactness suggests one to use a weak-type topology for $X$.

In conclusion, dealing with evolutionary problems,
we shall equip ${\cal G}$ with the evolutionary $\Gamma$-convergence
with respect to a weak-type topology for $X$.  

\bigskip
\noindent{\bf On weak $\Gamma$-convergence.\/}
Let the spaces $H, V, {\cal V}$ be as in \eqref{eq.BEN.spaces1} and \eqref{eq.BEN.spaces2}. 
For any $n$, let $\varphi_n\in {\cal E}_w(V)$ represent a semi-monotone operator 
$\alpha_n: V\to {\cal P}(V')$ as in \eqref{eq.SM1}---\eqref{eq.semim}. 
We define the measure $\mu$ be as in \eqref{eq.comstab.mu},  
the sequence $\{\widetilde\Phi_n\}$ in ${\cal V}_\mu \!\times\! L^2_\mu(0,T;V')$
as in \eqref{eq.comstab.BENfun'}, and study its (sequential) $\Gamma$-convergence. 

By confining $\{v_n^*\}$ to a bounded subset of $L^2_\mu(0,T;H)$, 
because of \eqref{eq.BEN.comp}
we are able to pass to the limit in $\int_0^T \langle v_n^*,v_n\rangle \, d\mu(t)$, as
$v_n\wto v$ in ${\cal V}_\mu$ and $v_n^*\wto v^*$ in $L^2(0,T;V')$,
up to extracting subsequences.
\footnote{ We denote the strong, weak, and weak star convergence respectively by $\to$,
$\wto$, $\wsto$. 
} 
Concerning the term $\int_0^T \! \langle D_tv,v\rangle \, d\mu(t)$, we have
\begin{equation}\label{eq.comstab.motiv3}
\int_0^T \! \langle D_tv,v\rangle \, d\mu(t) 
= {1\over2} \!\int_0^T \! D_t \big(\|v(t)\|_H^2\big) \, (T-t) \, dt 
= {1\over2} \!\int_0^T \!\! \|v\|_H^2 \, dt - {T\over2} \|v(0)\|_H^2,  
\end{equation} 
which by \eqref{eq.BEN.comp} is weakly continuous in ${\cal V}$. 
 
\begin{remark}\rm
At variance with $v\mapsto \int_0^T \! \langle D_tv,v\rangle \, d\mu(t)$, the functional 
$v\mapsto \int_0^T \! \langle D_tv,v\rangle \, dt = {1\over2} \|v(T)\|_H^2 -  \|v(0)\|_H^2$
is just weakly lower semicontinuous on ${\cal V}$, and the semicontinuity
is not sufficient for the (evolutionary) $\Gamma$-convergence of the sequence 
$\{\widetilde\Phi_n\}$. 
This is the main reason why we introduced the weight $\mu(t) = T-t$, 
which is tantamount to integrating a further time in $t$. 
\hfill$\Box$  
\end{remark}
  
\noindent{\bf The nonlinear weak topology.\/}
Dealing with the structural stability of the null-minimization of  
\eqref{eq.comstab.BENfun'}, one has to pass to the limit in the term
$\int_0^T \langle v_n^*- D_tv_n,v_n\rangle \, d\mu(t)$.
This induces us to complement the weak topology of 
$L^2_\mu(0,T;V \!\times\! V') \simeq L^2_\mu(0,T;V) \!\times\! L^2_\mu(0,T;V')$
with the convergence 
$\int_0^T \langle v_n^*,v_n \rangle \, d\mu(t) \to \int_0^T \langle v^*,v\rangle \, d\mu(t)$,
in view of the use of $\Gamma$-convergence.
This convergence will indeed provide the existence of a so-called recovery sequence.
Here it does not seem appropriate to use the product of 
the weak topology of $L^2_\mu(0,T;V)$ by the strong topology of $L^2_\mu(0,T;V')$.
Indeed, dealing with parabolic problems, by the BEN principle $u_n$ is typically paired 
with $u_n^* -D_tu_n$, and $\{D_tu_n\}$ converges just weakly in $L^2_\mu(0,T;V')$,
see Theorems~\ref{prop.repr} and \ref{prop.repr'}.

More specifically, let us still denote by $\pi$ the duality pairing of $V \!\times\! V'$:
$\pi(v,v^*) = \langle v^*,v\rangle$.
Along the lines of \cite{ViCalVar}, 
we shall name {\it nonlinear weak topology\/} of $V\!\times\! V'$, 
and denote by $\widetilde\pi$, the coarsest among the topologies of this space  
that are finer than the weak topology, and for which the mapping $\pi$ is continuous.
For any sequence $\{(v_n,v^*_n)\}$ in $V\!\times\! V'$, thus 
\begin{equation}
\begin{split}
&(v_n,v_n^*)\pito (v,v^*) \quad\hbox{in }V\!\times\! V'Ê
\quad\Leftrightarrow 
\\
&v_n\wto v \hbox{ \ \ in }V, \quad
v_n^*\wto v^* \hbox{ \ \  in }V', \quad
\langle v_n^*,v_n \rangle \to \langle v^*,v\rangle,
\end{split}
\end{equation}
and similarly for nets.  
This construction is extended to the space $L^2_\mu(0,T;V \!\times\! V')$ 
in an obvious way: in this case the duality product reads 
$L^2_\mu(0,T;V \!\times\! V')\to \erre: 
(v,v^*)\mapsto \int_0^T \langle v^*,v\rangle \, d\mu(t)$, and we set
\begin{equation}
\begin{split}
&(v_n,v_n^*)\pito (v,v^*) 
\quad\hbox{in }L^2_\mu(0,T;V \!\times\! V')Ê
\quad\Leftrightarrow 
\\
&v_n\wto v \hbox{ \ in }L^2_\mu(0,T;V), \;
v_n^*\wto v^* \hbox{ \  in }L^2_\mu(0,T;V') \; \hbox{ and }
\\
&\int_0^T \! \langle v_n^*,v_n \rangle \, d\mu(t) \to 
\! \int_0^T \! \langle v^*,v\rangle \, d\mu(t),
\end{split}
\end{equation}
and similarly for nets.

\bigskip
\noindent{\bf {\boldmath $\Gamma\widetilde\pi$}-compactness and 
{\boldmath $\Gamma\widetilde\pi$}-stability of ${\cal E}_{\widetilde\pi}(L^2_\mu(0,T;V))$.\/}
As the weak topology and the nonlinear weak topology $\widetilde\pi$ are nonmetrizable, 
some caution is needed in dealing with {\it sequential\/} 
$\Gamma$-convergence with respect to either topology;
for the former see e.g.\ \cite{At},\cite{Da}.
For functions defined on a topological space, the definition of 
$\Gamma$-convergence involves the filter of the neighborhoods of each point.
If the space is metrizable, that notion can equivalently be formulated in terms of the family of converging sequences; but this does not hold in general. 
We shall refer to these two notions as {\it topological\/} and {\it sequential\/} 
$\Gamma$-convergence, respectively.
Hereafter reference to the topological notion should be understood, if not otherwise stated.

Similarity between the nonlinear weak topology of $V\!\times\! V'$
and that of  $L^2_\mu(0,T;V \!\times\! V')$ is obvious.
Dealing with flows, here we are mainly concerned with the latter; 
for the sake of simplicity, we shall however develop our discussion for the former,
and leave the obvious reformulation to the reader.
 
It is known that bounded subsets of a separable and reflexive space equipped with the weak topology are metrizable.
The same holds for the nonlinear weak topology $\widetilde\pi$ of $V\!\times\! V'$, 
as it was proved in \cite{ViCalVar}. 
This property is at the basis of the next statement, where we define 
${\cal E}_{\widetilde\pi}(V)$ and ${\cal E}_{\widetilde\pi}(L^2_\mu(0,T;V))$
as in Section~\ref{sec.repres}, here with $\widetilde\pi$ instead of $\sigma$. 

\begin{theorem} ($\Gamma$-compactness and $\Gamma$-stability)  
\label{comsta} 
Let $V$ be a separable real Banach space, and $\{\gamma_n\}$ be an equi-coercive sequence in 
${\cal E}_{\widetilde\pi}(V)$, in the sense that 
\begin{equation}\label{eq.equicoer}
\sup_{n\in\enne} \big\{\|v\|_V + \|v^*\|_{V'}: 
(v,v^*)\in V\!\times\! V', \gamma_n(v,v^*) \le C \big\}<+\infty
\quad\forall C\in \erre.
\end{equation}
\indent 
Then:  
(i) There exists $\gamma: V\!\times\! V' \to \erre\cup\{+\infty\}$ such that, 
possibly extracting a subsequence, $\gamma_n$ $\Gamma\widetilde\pi$-converges 
to $\gamma$ both topologically and sequentially.
\hfill\break\indent
(ii) This entails that $\gamma \in {\cal E}_{\widetilde\pi}(V)$. 
Moreover $\gamma \in {\cal F}(V)$ if $\gamma_n \in {\cal F}(V)$ for any $n$.
\hfill\break\indent
(iii) If $\alpha_n$ ($\alpha$, resp.) is the operator that is represented by 
$\gamma_n$ ($\gamma$, resp.) for any $n$, then 
\begin{equation}\label{eq.comstab.lim}
\forall \text{ sequence }\{(v_n,v_n^*)\in \gr(\alpha_n)\}, \quad
(v_n,v^*_n)\pito (v,v^*)
\quad\Rightarrow\quad v^*\in \alpha(v).
\end{equation}
That is, the superior limit in the sense of Kuratowski of the graph of the $\alpha_n$s is 
included in the graph of $\alpha$.
\end{theorem} 

\noindent{\bf Proof.\/}  
Part (i) is Theorem~4.4 of \cite{ViCalVar}.
 
As the functional $\pi$ is obviously $\widetilde\pi$-continuous,
the sequential $\widetilde\pi$-lower semicontinuity of the $\gamma_n$s and the property 
``$\gamma_n \ge \pi$'' are preserved by $\Gamma\widetilde\pi$-convergence;
thus $\gamma \in {\cal E}_{\widetilde\pi}(V)$. 
If $\gamma_n \in {\cal F}(V)$ for any $n$, namely if the $\gamma_n$s are also convex,
then the same holds for $\gamma$.
Part (ii) is thus established. 

Next let the operators $\{\alpha_n\}$ and $\alpha$ be as prescribed in part (iii), 
$(v_n,v^*_n)\in \gr(\alpha_n)$ for any $n$, and $(v_n,v^*_n)\pito (v,v^*)$. 
By \eqref{eq.fitzp.nonconv3}, thus $ \gamma_n(v_n, v^*_n) = \langle v^*_n,v_n \rangle$ 
for any $n$. Recalling the definition of $\Gamma\widetilde\pi$-convergence, 
if $(v_n,v^*_n)\pito (v,v^*)$ we then have
\begin{equation}\label{eq.comstab.lim'}
\gamma(v,v^*) 
\le\liminf_{n\to \infty} \gamma_n(v_n, v^*_n) = 
\liminf_{n\to \infty} \; \langle v^*_n,v_n \rangle = \langle v^*,v \rangle.
\end{equation}
Thus $v^*\in \alpha(v)$, as $\gamma$ represents $\alpha$. 
The implication \eqref{eq.comstab.lim} is thus established.   
\hfill$\Box$
 
\begin{remarks}\label{eq.comstab.rem}\rm 
(i) By the same argument, Theorem~\ref{comsta} holds also if the space $V$ and 
$\gamma_n\in {\cal E}_{\widetilde\pi}(V)$ are respectively replaced by 
$L^2_\mu(0,T;V)$ and $\gamma_n \in {\cal E}_{\widetilde\pi}(L^2_\mu(0,T;V))$.
\footnote{ This is the case of interest for the next two sections.
}

(ii) In general the sequence of operators $\{\alpha_n\}$ does not converge in the sense of 
Kuratowski. 
For instance, let us define $F_b$ as in \eqref{eq.fitzp.counterex} for any $b\ge 1/2$.
If $b_n> 1/2$ for any $n$ and $b_n\to 1/2$, then $F_{b_n}$ 
$\Gamma\widetilde\pi$-converges to $F_{1/2}$, 
but the represented operators $\alpha_n$ do not converge in the sense of Kuratowski
to the operator $\alpha$ that is represented by $F_{1/2}$.
Actually, in this case $\gr(\alpha_n) = \{(0,0)\}$, so that
the inferior limit in the sense of Kuratowski of the graph of the $\alpha_n$s does 
not include the graph of $\alpha$.  
\end{remarks}  
  
\noindent{\bf Onset of long memory in the limit?} 
Let $\{\varphi_n\}$ be a sequence of representative functions of 
${\cal E}_{\widetilde\pi}(V)$, and define the superposition 
(i.e., Nemytski\u{\i}-type) operators 
\footnote{ These operators should not be mixed up with the functions 
$\gamma_n\in {\cal E}_{\widetilde\pi}(L^2_\mu(0,T;V))$ of the 
Remark!\ref{eq.comstab.rem} (i).
}
\begin{equation}
\psi_n: L^2_\mu(0,T;V \!\times\! V')\to L^1(0,T): w\mapsto \varphi_n(w(\cdot)) \quad\forall n.
\end{equation} 
The following question arises:
\begin{equation}
\begin{split}
&\text{ if $\psi_n$ $\Gamma\widetilde\pi$-converges to some operator $\psi$ in the sense of 
\eqref{eq.evol.defgamma.1+}, }
\\
&\text{ is then $\psi$ necessarily a superposition operator, too? }
\end{split}
\end{equation}  
In other terms, there exists then a mapping $\varphi:  V \!\times\! V'\to \erre \cup\{+\infty\}$ 
such that $\psi_w = \varphi(w(\cdot))$ for any $w\in L^2_\mu(0,T;V \!\times\! V')$?
This would exclude the onset of long memory in the limit,
as it was discussed in \cite{ViCalVar}.
A positive answer is provided by the next statement, 
which is essentially a particular case of Theorem~\ref{teo.comp}, 
augmented with the hypothesis \eqref{eq.comstab.reppn} and with the corresponding property 
in the thesis. 
Here we just deal with the Hilbert set-up and with $p=2$ for technical reasons. 
We allow the representing functions $\varphi_n$ to depend explicitly on time, too,
since this does not raise any further difficulty. 

\begin{theorem} \label{teo.comp'}
Let $V$ be a real separable Hilbert space,
and $\mu$ be the measure on $]0,T[$ such that $d\mu(t) = (T-t) \, dt$.
Let $\{\varphi_n\}$ be a sequence of normal functions 
$]0,T[{} \!\times\! V \!\times\! V'\to \erre^+$ 
(i.e., for any $n$, $\varphi_n$ is globally measurable, 
and $\varphi_n(t,\cdot)$ is lower semicontinuous for a.e.\ $t\in {}]0,T[$), such that   
\begin{eqnarray} 
&\varphi_n(t,\cdot) \in {\cal E}_{\widetilde\pi}(V)
\qquad\hbox{for a.e.\ }t\in {}]0,T[, \forall n,  
\label{eq.comstab.reppn}
\\
&\begin{split}
&\exists C_1,C_2,C_3 >0: \forall n,\hbox{for a.e.\ }t\in {}]0,T[, \forall w\in V \!\times\! V', 
\\
&C_1 \|w\|_{V \!\times\! V'}^2\le \varphi_n(t,w) \le C_2\|w\|_{V \!\times\! V'}^2 +C_3, 
\end{split}   
\label{eq.comstab.equibc'}  
\\
&\varphi_n(t,0) =0 \qquad\hbox{ for a.e.\ }t\in {}]0,T[,\forall n,
\label{eq.comstab.nul'}
\end{eqnarray} 
and define the operators $\psi_n: L^2_\mu(0,T;V \!\times\! V')\to L^1_\mu(0,T)$ by
\begin{equation}\label{eq.comstab.super'} 
\psi_{n,w}(t) = \varphi_n(t,w(t))  
\qquad\forall w\in L^2_\mu(0,T;V \!\times\! V'),\hbox{ for a.e.\ }t\in {}]0,T[,\forall n.
\end{equation}  
\indent
Then there exists a normal function $\varphi: {}]0,T[{} \times V \!\times\! V'\to \erre^+$ 
such that
\begin{equation}\label{eq.comstab.repp}
\begin{split}
\varphi(t,\cdot) \in {\cal E}_{\widetilde\pi}(V)
\quad
(\varphi(t,\cdot) \in {\cal F}(V) \text{ if }
\varphi_n(t,\cdot) \in {\cal F}(V) \text{ for any }n)
\\
\hbox{for a.e.\ }t\in {}]0,T[,
\end{split}
\end{equation}  
and such that, defining the corresponding operator 
$\psi: L^2_\mu(0,T;V \!\times\! V')\to L^1_\mu(0,T)$ 
as in \eqref{eq.comstab.super'}, possibly extracting a subsequence,
\begin{equation}\label{eq.comstab.tesi'} 
\begin{split}
&\hbox{ $\psi_n$ sequentially $\Gamma$-converges to }\psi
\\
&\hbox{ in the topology $\widetilde\pi$ of $L^2_\mu(0,T;V \!\times\! V')$ and } 
\\
&\hbox{ in the weak topology of $L^1_\mu(0,T)$
(cf.\ \eqref{eq.evol.defgamma.1+}). }
\end{split}
\end{equation}
\indent
Finally, if $\varphi_n$ does not depend on $t$ for any $n$, 
then the same holds for $\varphi$. 
\end{theorem}
 
\noindent{\bf Proof.\/}
Let us apply Theorem~\ref{teo.comp} with $X= V \!\times\! V'$,
$p=2$ and the topology $\tau = \widetilde\pi$; the hypothesis \eqref{eq.evol.gamcom} 
is indeed fulfilled for this topology, because of part (i) of Theorem~\ref{comsta}. 
It then suffices to show that \eqref{eq.comstab.reppn} entails the existence of 
a normal representative function $\varphi$, as it is stated in \eqref{eq.comstab.repp}.
In order to prove that $\varphi$ is a representative function 
(see \eqref{eq.fitzp.nonconv=}), let us set
\begin{equation}\label{eq.comstab.nonconv=}
\begin{split}
&J_n(t,v,v^*) := \varphi_n(t,v,v^*) - \langle v^*,v\rangle
\\[1mm]
&J(t,v,v^*) := \varphi(t,v,v^*) - \langle v^*,v\rangle 
\end{split}
\qquad\forall (v,v^*)\in V \!\times\! V', \hbox{ for a.e.\ }t,\forall n.
\end{equation}  
By \eqref{eq.comstab.reppn}, 
\begin{equation}\label{eq.comstab.nonconv+}
\int_A J_n(t,v,v^*) \, d\mu(t) \ge 0 
\qquad\forall (v,v^*)\in V \!\times\! V', \forall A\in {\cal L}(0,T),
\end{equation} 
and by \eqref{eq.comstab.tesi'} this inequality is preserved in the limit. 
Therefore $J(t,v,v^*)\ge 0$ a.e.\ in $]0,T[$.
As $\varphi(t,\cdot)$ is $\widetilde\pi$-lower semicontinuous, 
\eqref{eq.comstab.repp} follows.

The final statement about $t$-independence stems from the analogous assertion of 
Theorem~\ref{teo.comp}.
\hfill$\Box$

\section{Structural properties of semi-monotone flows}
\label{sec.par}

\noindent
In this section we prove the structural compactness and stability of flows of the form 
$D_tu + \alpha(u) \ni h$ for semi-monotone operators $\alpha$. 
 
\bigskip
\noindent{\bf Quasilinear parabolic operators in abstract form.\/}  
Let the Hilbert spaces $V$, $H$, ${\cal V}$ be as in Section~\ref{sec.BEN},
and three sequences $\{b_n\},\{\beta_n\},\{f_n\}$ be given such that 
\begin{eqnarray} 
&b_n: V\to {\cal P}(V')\text{ is maximal monotone, }\forall n,  
\label{eq.par.op0} 
\\
&\beta_{n,z}: V\to V'\text{ is maximal monotone, }\forall z\in V,\forall n, 
\label{eq.par.op1} 
\\
&V_w\to (V')_s: z\mapsto \beta_{n,z}(v)\text{ is continuous, }\forall v\in V,\forall n, 
\label{eq.par.op2}
\\ 
&g_{n,z} \text{ is the Fitzpatrick function of }b+ \beta_{n,z}, \forall z\in V, \forall n.
\label{eq.par.op3} 
\end{eqnarray} 
Setting
\begin{eqnarray}\label{eq.par.opm} 
&\alpha_n(v) := b_n(v) + \beta_{n,v}(v) \qquad\forall v\in V,\forall n, 
\label{eq.par.opm1} 
\\
&\varphi_n(v,v^*) := g_{n,v}(v,v^*)
\qquad\forall (v,v^*)\in V\times V',\forall n, 
\label{eq.par.opm2} 
\end{eqnarray}
by Theorem~\ref{BEN.param} then $\varphi_n\in {\cal E}_{\widetilde\pi}(V)$ and it 
$\widetilde\pi$-represents $\alpha_n$ for any $n$.  

We shall assume that
\begin{eqnarray}  
&\exists C_1, C_2>0: \forall n, 
\forall v\in V, \quad 
\langle \alpha_n(v),v \rangle \ge C_1|v\|_V^2 - C_2, 
\label{eq.par.alpha2}
\\ 
&\exists C_3, C_4>0: \forall n, 
\forall v\in V, \quad 
\|\alpha_n(v)\|_{V'} \le C_3 \|v\|_V + C_4,
\label{eq.par.alpha3} 
\\
&0\in \alpha_n(0) \qquad\forall n. 
\label{eq.par.alpha3'}
\end{eqnarray} 
The condition \eqref{eq.par.alpha3'} is not really restrictive: 
if it is not satisfied, it can be recovered by selecting any $a\in\alpha_n(0)$ and 
then replacing $\alpha_n$ by $\bar\alpha_n = \alpha_n - a$. 
Let two sequences $\{u^0_n\}$ and $\{h_n\}$ be also given such that  
\begin{eqnarray} 
&u^0_n\to u^0 \quad\hbox{ in }H,
\label{eq.par.alpha4}
\\
&h_n\to h \quad\hbox{ in }L^2(0,T;V').
\label{eq.par.alpha5}
\end{eqnarray} 
 
We are now able to introduce the following initial-value problem, for any $n$, 
\begin{equation}\label{pbn1}  
\left\{\begin{split} 
&u_n\in {\cal V},
\\
&D_tu_n + \alpha_n(u_n) \ni h_n \qquad \hbox{ in $V'$, a.e.\ in }{}]0,T[,
\\
&u_n(0) = u^0_n.
\end{split}\right.
\end{equation}

\begin{lemma} [\cite{BrHe}] \label{regul}
Under the hypotheses above, for any $n$ the initial-value problem~\eqref{pbn1} has one and only one solution $u_n\in {\cal V}$. Moreover, the sequence $\{u_n\}$ is bounded in 
${\cal V}$.
\end{lemma} 
  
Let us define the measure $\mu$ as in \eqref{eq.comstab.mu}. 

\begin{theorem}(Structural compactness and stability)\label{eq.par.stru}
Let \eqref{eq.par.op0}--\eqref{eq.par.alpha5} be fulfilled, 
and for any $n$ let $u_n$ be a solution of the Cauchy problem~\eqref{pbn1}. Then:
\\
\indent
(i) There exists $u\in {\cal V}$ such that, possibly extracting a subsequence,  
\begin{equation}\label{eq.par.1a}
u_n\wto u \qquad\hbox{ in }{\cal V}.
\end{equation} 
\indent
(ii) There exists a function $\varphi\in {\cal E}_{\widetilde\pi}(V)$ such that, setting 
\begin{eqnarray}   
&\psi_{n,w}(t) = \varphi_n(w(t)), \qquad
\psi_w(t) = \varphi(w(t)),  
\\[1mm]
&\hbox{for a.e.\ }t\in{}]0,T[,\forall w\in L^2_\mu(0,T;V \!\times\! V'), \forall n,
\label{eq.par.de''}
\end{eqnarray}
then $\psi_n,\psi: L^2_\mu(0,T;V \!\times\! V') \to L^1_\mu(0,T)$, and, 
possibly extracting a subsequence, 
\begin{equation}\label{eq.par.tesi+}
\begin{split}
&\hbox{ $\psi_n$ sequentially $\Gamma$-converges to }\psi
\\
&\hbox{ in the topology $\widetilde\pi$ of $L^2_\mu(0,T;V \!\times\! V')$ and } 
\\
&\hbox{ in the weak topology of $L^1_\mu(0,T)$ 
(cf.\ \eqref{eq.evol.defgamma.1+}).}
\end{split}
\end{equation} 
\indent
(iii) Denoting by $\alpha: V\to {\cal P}(V')$ the operator that is represented by $\varphi$,
$u$ solves the corresponding initial-value problem
\begin{equation}\label{pb=}  
\left\{\begin{split} 
&u\in {\cal V},
\\
&D_tu + \alpha(u) \ni h \qquad \hbox{ in $V'$, a.e.\ in }{}]0,T[,
\\
&u(0) = u^0.
\end{split}\right.
\end{equation}
\end{theorem}

\medskip
\noindent{\bf Proof.\/} 
For any $n$, let us define the affine subspace ${\cal V}_{\mu,u_n^0}$ of ${\cal V}$
as in \eqref{eq.comstab.spaces3}, and the sequence of nonnegative twice-time-integrated functionals $\{\widetilde\Phi_n\}$ as in \eqref{eq.comstab.BENfun'}. 

For any $n$, by part (ii) of Theorem~\ref{prop.repr'}, $u_n$ solves 
the initial-value problem~\eqref{pbn1} if and only if  
\begin{equation}\label{eq.par.min1}  
u_n\in {\cal V}_{\mu,u_n^0}, 
\qquad
\widetilde\Phi_n(u_n,h_n) = 0\;\; 
\big(\! =\inf_{{\cal V}_{\mu,u_n^0}} \widetilde\Phi_n(\cdot,h_n) \big).
\end{equation} 

By Lemma~\ref{regul}, \eqref{eq.par.1a} holds up to extracting a subsequence.
By \eqref{eq.comstab.motiv3} then
\begin{equation}\label{eq.par.co1}
\int_0^T \langle h_n - D_tu_n, u_n\rangle \, d\mu(t) \to
\int_0^T \langle h - D_tu, u\rangle \, d\mu(t),
\end{equation}
so that
\begin{equation}\label{eq.par.co2}
(h_n - D_tu_n, u_n)\pito (h - D_tu, u)
\qquad\hbox{ in }L^2_\mu(0,T;V \!\times\! V').
\end{equation}

Because of \eqref{eq.par.alpha2} and \eqref{eq.par.alpha3},
we can apply Theorem~\ref{teo.comp'}.
As here the functions $\varphi_n$s do not explicitly depend on time, 
there exists a function $\varphi\in {\cal E}_{\widetilde\pi}(V)$ such that,
defining $\psi$ as in \eqref{eq.par.de''}, 
$\psi: L^2_\mu(0,T;V \!\times\! V') \to L^1_\mu(0,T)$ and, 
possibly extracting a subsequence, \eqref{eq.par.tesi+} is fulfilled. 
This yields
\begin{equation}\label{eq.par.ga} 
\int_0^T \! \varphi_n(v,v^*) \, d\mu(t) 
\Gammatopi \int_0^T \! \varphi(v,v^*) \, d\mu(t)
\qquad\hbox{ sequentially in }V \!\times\! V'.
\end{equation}
 
Therefore
\begin{equation}
\begin{split}
&\int_0^T \! \varphi(u, h- D_tu) \, d\mu(t)  
\overset{\eqref{eq.par.co2},\eqref{eq.par.ga}}{\le}
\liminf_{n\to \infty} \int_0^T \! \varphi_n(u_n, h_n- D_tu_n) \, d\mu(t)
\\
&\overset{\eqref{eq.par.min1}}{=}
\!\! \liminf_{n\to \infty} \; \int_0^T \langle h_n - D_tu_n, u_n\rangle \, d\mu(t) 
\overset{\eqref{eq.par.co1}}{=}
\!\int_0^T \! \langle h - D_tu, u\rangle \, d\mu(t) .  
\end{split}
\end{equation}
Hence, defining $\widetilde\Phi$ as in \eqref{eq.comstab.BENfun'}, 
\begin{equation}\label{eq.par.min2}  
u\in {\cal V}_{\mu,u^0}, 
\qquad
\widetilde\Phi(u,h) = 0\;\; 
\Big(\! =\inf_{{\cal V}_{\mu,u^0}} \widetilde\Phi(\cdot,h) \Big).
\end{equation}

By Theorem~\ref{prop.repr'}, \eqref{pb=} is thus established. 
\hfill$\Box$ 

\begin{remarks}\rm
(i) This theorem extends the results of Section~8 of \cite{ViCalVar} to semi-monotone flows, 
and improves them since it expresses the limit functional in terms of the superposition function 
$\varphi$. This excludes the onset of long memory.

(ii) We just saw that the variational formulation of the flow of semi-monotone operators is preserved in the $\Gamma$-limit. 
On the other hand, if the $\alpha_n$s are cyclical maximal operators, 
the structure of gradient flow need not be preserved by $\Gamma$-convergence.

More specifically, let us assume that $\alpha_n = \partial \varphi_n$, where $\varphi_n$ is a
lower semicontinuous convex function $V \!\times\! V'\to \erre\cup \{+\infty\}$ for any $n$. 
The cyclical maximal monotone operator $\alpha_n$ is thus represented by the Fenchel function $f_n: (v,v^*)\mapsto \varphi_n(v)+\varphi_n^*(v^*)$.  
But $f_n\Gammatopi f$ does not entail that $f$ be a Fenchel function, at variance with what happens if  $f_n\Gammato f$ either in $V_s \!\times\! V'_w$ or in $V_w \!\times\! V'_s$.
This is illustrated in Section~5 of \cite{ViCalVar}, where a counterexample is also displayed. 
  
(iii) Theorem~\ref{eq.par.stru} can be extended in several ways. 
For instance, a time-dependent semi-monotone operator $\alpha: V\!\times {}]0,T[{} \to 
{\cal P}(V')$ may be represented by a time-dependent representative function 
$f_\alpha(\cdot,\cdot,t)\in {\cal E}_{\widetilde\pi}(V)$.
An existence result analogous to Lemma~\ref{regul} holds for nonautonomous
semi-monotone operators, and Theorem~\ref{eq.par.stru} takes over to this more general set-up.
In this case $\psi$ explicitly depends on time;
so the generality of Theorems~\ref{teo.comp'} is fully exploited. 

(iv) If instead of prescribing an initial condition one assumes the solution $u$
to be $T$-periodic in time, then $\int_0^T \langle D_tu,u\rangle \, dt =0$.
In this case it is not necessary to introduce any weight function $\mu$,
and the analysis is much simpler; see \cite{ViCalVar}.
\hfill$\Box$  
\end{remarks}

\section{Application to PDEs of mathematical physics} 
\label{sec.PDE}   

\noindent 
Theorem~\ref{eq.par.stru} provides the structural compactness 
and stability of the Cauchy problem associated with several quasilinear PDEs. 
In this section we illustrate some examples issued from mathematical physics. 

\medskip
(i) Quasilinear diffusion in a bounded Lipschitz domain $\Omega$ of $\erre^N$ ($N\ge 1$)
can be represented by the quasilinear parabolic equation 
\begin{equation}\label{eq.spaces.1}
D_tu - \nabla \!\cdot\! g(x,u,\nabla u) \ni h
\qquad\hbox{ in $Q:= \Omega \!\times\! {}]0,T[$ } \; (\nabla \cdot\! := \div).
\end{equation}  

More precisely, let us assume that $g(x,v,\xi) = g_1(\xi) + g_2(x,v,\xi)$,
\begin{equation}\label{eq.spaces.1=}
\begin{split}
&g_1: \erre^N\to {\cal P}(\erre^N)
\qquad
g_2: \Omega \!\times\! \erre \!\times\! \erre^N\to \erre^N, 
\\
&v\mapsto g_1(v)\text{ \ is maximal monotone, }
\\ 
&x\mapsto g_2(x,v, \xi)\text{ \ is Lebesgue-measurable, }\forall (v,\xi), 
\\
&v\mapsto g_2(x,v, \xi)\text{ \ is continuous, for a.e.\ }x, \forall \xi,
\\
&\xi \mapsto g_2(x,v, \xi)\text{ \ is monotone, for a.e.\ }x, \forall v,
\\
&\begin{split}
&\exists c_1,c_2>0: \text{for a.e.\ }x, \forall (v,\xi), \;\;\;
\\
&|g_1(\xi)| + |g_2(x,v, \xi)| \le c_1 (|v| + |\xi|) + c_2,
\end{split}
\end{split}
\end{equation} 
and set 
\begin{equation}\label{eq.spaces.1b}
\alpha(v) = -\nabla \cdot [g_1(v)] -\nabla \cdot [g_2(x,v,\nabla v)]
\qquad\text{ in }{\cal D}'(\Omega),\forall v\in H^1_0(\Omega).
\end{equation}
Then $\alpha:H^1_0(\Omega) \to {\cal P}(H^{-1}(\Omega))$ is semi-monotone. 
In this case, prescribing e.g.\ the homogeneous Dirichlet condition, under standard assumptions
of coerciveness the existence of a weak solution of the flow associated with \eqref{eq.spaces.1}
follows from results of \cite{Br0},\cite{BrHe},\cite{Li}.
Setting $H= L^2(\Omega)$ and $V= H^1_0(\Omega)$, Theorem~\ref{eq.par.stru} 
then provides structural compactness and stability.

\medskip
(ii) If $N=3$ one can also deal with the quasilinear parabolic vector equation
\begin{equation}\label{eq.spaces.2}
\mu D_tH + \nabla \!\times\! g(x, H, \nabla \!\times\! H) \ni 0
\qquad\hbox{ in $Q$ }\; (\nabla \times\! := \curl). 
\end{equation}
This equation arises e.g.\ by coupling the Faraday law of magnetic induction, 
$D_tB + \nabla \!\times\! E = 0$, 
with a nonlinear Ohm law $E \ni g(x,H,J)$, 
with the constitutive relation $B = \mu H$,
and with the Amp\`ere law $J = \nabla \!\times\! H$.
The latter is here written neglecting the displacement current, 
by the so-called {\it eddy current approximation.\/} 
The dependence of $E$ on $H$ accounts for the Hall effect. 

In this case, denoting the outward-oriented unit normal vector-field on $\partial\Omega$ 
by $\nu$, we assume properties analogous to \eqref{eq.spaces.1=} and set
\begin{equation}\label{eq.spaces.3}
\begin{split}
&H= \big\{v\in L^2(\Omega)^3: \nabla \!\cdot\! v =0 
\hbox{ in }{\cal D}'(\Omega) \big\},
\\
&V= \big\{v\in H: \nabla \!\times\! v \in L^2(\Omega)^3, \;
\nu \!\times\! v = 0 \hbox{ in } H^{-1/2}(\partial\Omega)^3 \big\},  
\\
&\alpha(v) = \nabla \!\times\! g_1(\nabla \!\times\! v) 
+ \nabla \!\times\! g_2(x,v,\nabla \!\times\! v) 
\qquad\forall v\in V.
\end{split}
\end{equation}  
See e.g.\ Chap.\ IV of \cite{Vi96} and references therein. 
The operator $\alpha: V\to {\cal P}(V')$ is semi-monotone, and 
structural compactness and stability are provided by Theorem~\ref{eq.par.stru}. 

\medskip
(iii) If $k$ is a positive function of Borel class, then the equation 
\begin{equation}\label{eq.heat1}
D_t u- \nabla\cdot [k(u) \nabla u] =h,
\qquad\text{ in }Q
\end{equation}
with $k$ positive and integrable, is a particular case of \eqref{eq.spaces.1}, but can also be addressed via a different functional set-up.

Defining the {\it Kirchhoff transform\/} $\theta(v) = \int_0^v k(s) \, ds$ for any $v\in\erre$,
the equation \eqref{eq.heat1} reads
\begin{equation}\label{eq.heat2}
D_t u- \Delta\theta(u) =h
\qquad\text{ in }Q.
\end{equation} 
Although the function $\theta$ is maximal monotone, it is promptly checked that 
in general the operator 
$H^1_0(\Omega)\to H^{-1}(\Omega): v\mapsto - \Delta\theta(v)$ is not even monotone.
However, by a classical construction known as {\it change of pivot space,\/} 
see e.g.\ \cite {Li} p.~190 and \cite{Vi96} Section~II.6, one can reduce the equation 
\eqref{eq.heat2} to the form $D_tu + \alpha(u)= h$ with $\alpha$ maximal monotone.

\medskip
(iv) The weak formulation of the two-phase {\it Stefan problem\/}
(the classical mathematical model of phase transitions) is also of the form \eqref{eq.heat2}, 
for a maximal-monotone operator $\theta$; see e.g.\ \cite{Vi96}.
In this case $u$ and $\theta =\theta(u)$ respectively represent the density of internal energy 
(of the enthalpy, if pressure is maintained constant) and the temperature.
This can easily be extended to account for nonlinear diffusion. 

\medskip
(v) In presence of a prescribed velocity field $v$ ($\in L^2(\Omega)^3$, say), 
diffusion with convection can be represented by inserting the convective term 
$v \!\cdot\! \nabla u = \{\sum_{j=1}^N v_j D_{x_j} u\}$ into \eqref{eq.heat1}:
\begin{equation}\label{eq.heatconv1}
D_t u + v \!\cdot\! \nabla u- \nabla\cdot [k(u) \nabla u] =h.
\qquad\text{ in }Q
\end{equation}
If $\nabla \!\cdot\! v \le 0 \hbox{ in }{\cal D}'(\Omega)$, 
i.e.\ the medium is not expanding, then the linear operator 
$H^1_0(\Omega)\to H^{-1}(\Omega): u\mapsto v \!\cdot\! \nabla u$ \ is
monotone, since 
\begin{equation}\label{eq.heatconv2}
\int_\Omega (v \cdot \nabla u) u \, dx 
= {1\over2} \sum_{j=1}^N \int_\Omega v_j D_{x_j} u^2 \, dx
= - {1\over2} \sum_{j=1}^N \int_\Omega (D_{x_j}v_j) u^2 \, dx \ge 0
\end{equation}
for any $u\in H^1_0(\Omega) - H^1_0(\Omega)$.
The operator $u\mapsto v \cdot \nabla u- \nabla\cdot [k(u) \nabla u]$
is thus semi-monotone. 

\medskip
(vi) Next we outline a variant of the Navier-Stokes model for the
flow of an incompressible non-Newtonian viscous fluid;
see e.g.\ \cite{DuLi},\cite{MaRa},\cite{Te} and references therein.
Let $\Omega$ be a bounded domain of $\erre^N$ ($N=2$ or $3$).
Let us denote the velocity field by $u$, and its symmetrized gradient by $\nabla^su$.
Let us assume that the viscous stress $\sigma$ and the strain-rate are related as follows:
\begin{equation}\label{eq.NS1}
\sigma \in \gamma(\nabla^su),
\quad\text{ with }\quad
\gamma: \erre^{N \!\times\! N} \to {\cal P}(\erre^{N \!\times\! N})\text{ maximal monotone. }
\end{equation} 
We denote the pressure field by $p$, the constant density by $\rho$, 
a prescribed force field by $h$, and consider the system
\begin{equation}\label{eq.NS2}
\left\{\begin{split}
&\rho D_tu + \rho u \cdot\! \nabla u - \nabla \!\cdot\! \gamma(\nabla^su) 
\ni h -\nabla p \qquad\hbox{ in }Q,
\\[1mm]
&\nabla \!\cdot\! u =0 \qquad\quad\hbox{ \ in }Q,
\\[1mm]
&u(\cdot,0) = u\,{}^0 \qquad\hbox{ in }\Omega.
\end{split}
\right.
\end{equation}

Let us prescribe e.g.\ the homogeneous Dirichlet condition, set 
\begin{equation}\label{eq.NS3}
H= \big\{v\in L^2(\Omega)^N: \nabla \!\cdot\! v =0 \hbox{ in }{\cal D}'(\Omega)\big\},
\qquad
V= H\cap H^1_0(\Omega)^N,
\end{equation} 
and define the maximal monotone operator ${\widetilde\gamma}$ associated to the mapping 
$\gamma$:
\begin{equation}
{\widetilde\gamma}: V\to {\cal P}(V'): v\mapsto - \nabla \!\cdot\! \gamma(\nabla^su).
\end{equation}

Existence of a weak solution of the flow \eqref{eq.NS2} can then be proved via the standard procedure.   
By the Fitzpatrick theorem, if $f_{\widetilde\gamma}$ is defined as in \eqref{eq.BEN.=},  
then the flow \eqref{eq.NS2} is equivalent to the null-minimization of the functional
\begin{equation}\label{eq.NS5} 
\begin{split}
&\widetilde\Psi(v,v^*) := \!\! \int_0^T 
\!\! \big[f_{\widetilde\gamma}(v,h - \rho v \cdot\! \nabla v- \rho D_tv) 
- \langle h,v\rangle\big] \, d\mu(t)  
\\
&\qquad\qquad
+{\rho\over2}\! \int_0^T \!\!\! \|v(t)\|_H^2 \, dt - {\rho T\over2} \|u\,{}^0\|_H^2 
\qquad
\forall (v,v^*)\in {\cal V}_{\mu,u\,{}^0} \!\times\! L^2(0,T;V'), 
\\
&\widetilde\Phi(v,v^*) = +\infty 
\qquad\qquad\qquad\qquad\qquad\qquad\;
\hbox{ for any other }(v,v^*)\in {\cal V} \!\times\! {\cal V}'.
\end{split} 
\end{equation} 

Theorem~\ref{eq.par.stru} provides the structural compactness and stability of this model.   

\medskip
(vii) Analogous conclusions can be attained for the flow associated to quasilinear hyperbolic equations, like e.g.\
\begin{equation}\label{eq.iperb}  
D_t^2u - \Delta u - D_t \nabla \!\cdot\! g(u,\nabla u) = h
\qquad\hbox{ in }Q,
\end{equation}
with $L^2(\Omega)^N \!\times\! H^1_0(\Omega)^N\to H^{-1}(\Omega)^3: 
(u,v)\mapsto - \nabla \!\cdot\! g(v,\nabla u)$ 
continuous with respect to $v$ and maximal monotone with respect to $u$.
More generally, $g$ can be replaced by a multi-valued mapping of the form 
\eqref{eq.spaces.1=}.

By integrating \eqref{eq.iperb} in time, one gets the integro-differential equation
\begin{equation}  
\begin{split}  
&D_tu - \Delta \!\!\int_0^t \! u(x,s) \, ds -\nabla \!\cdot\! g(u,\nabla u) 
\\
&= \! \int_0^t\! h(x,s) \, ds + D_tu(\cdot,0) - \nabla \!\cdot\! g(u^0,\nabla u^0)
\qquad\hbox{ in }Q.
\end{split}
\end{equation} 
One can then exploit the semi-monotonicity of the operator 
$u\mapsto - \nabla \!\cdot\! g(u,\nabla u)$.  

\bigskip
In conclusion, by Theorem~\ref{eq.par.stru}, the Cauchy problem
associated with any of the equations of this section is structurally compact and structurally stable.
This requires minimal regularity hypotheses on the data, essentially
the same assumptions that provide existence of a weak solution.

\section{Appendix. Evolutionary {\boldmath $\Gamma$}-convergence of weak type}
\label{sec.evol}  

\noindent
In this section we extend De Giorgi's notion of $\Gamma$-convergence to operators 
(rather than functionals) that act on time-dependent functions with range in a Banach space 
$X$, and state a result of $\Gamma$-compactness, along the lines of \cite{Vi17}. 
This is instrumental to Theorem~\ref{teo.comp'} of Section~\ref{sec.comstab}. 

\bigskip
\noindent{\bf Functional set-up.}
Let $X$ be a real separable and reflexive Banach space, $p\in [1,+\infty[$, $T>0$, and 
define the measure $\mu$ as in \eqref{eq.comstab.mu}. 
Let us equip $L^p_\mu(0,T;X)$ with a topology $\tau$ that is finer than the weak topology.
\footnote{ We assume this because of the application of Sections~\ref{sec.par}.
A reader interested just in evolutionary $\Gamma$-convergence might go through this section assuming that $\mu$ is the Lebesgue measure and that $\tau$ is the weak topology.
}  
For any operator $\psi: L^p_\mu(0,T;X)\to L^1_\mu(0,T):w\mapsto \psi_w$, let us set 
\begin{equation}\label{eq.evol.crochet}
[\psi,\xi](w) = \int_0^T \psi_w(t) \, \xi(t) \, d\mu(t) 
\qquad\forall w\in L^p_\mu(0,T;X), \forall \xi\in L^\infty(0,T). 
\end{equation}  
 
\noindent{\bf Evolutionary $\Gamma$-convergence of weak type.}
Let $\{\psi_n\}$ be a sequence of operators $L^p_\mu(0,T;X)\to L^1_\mu(0,T)$,
such that for any bounded subset $A$ of $L^p_\mu(0,T;X)$
the set $\{\psi_{n,w}: w\in A, n\in \enne\}$ is bounded in $L^1_\mu(0,T)$. 
If $\psi$ also is an operator $L^p_\mu(0,T;X)\to L^1_\mu(0,T)$, we shall say that 
\begin{equation}\label{eq.evol.defgamma.1+}
\begin{split}
&\hbox{ $\psi_n$ sequentially $\Gamma$-converges to }\psi
\\
&\hbox{ in the topology $\tau$ of $L^p_\mu(0,T;X)$ and } 
\\
&\hbox{ in the weak topology of }L^1_\mu(0,T)
\end{split}
\end{equation}
if and only if 
\begin{equation}\label{eq.evol.defgamma.2=}
\begin{split}
&[\psi_n, \xi] \text{ sequentially $\Gamma\tau$-converges to $[\psi,\xi]$ in }
\\[2mm]
&L^p_\mu(0,T;X), \; \forall \text{ nonnegative } \xi\in L^\infty(0,T).
\end{split} 
 \end{equation} 
 
\begin{remark}\rm
The present definition of {\it evolutionary $\Gamma$-convergence of weak type\/} 
is not equivalent 
either to that of \cite{SaSe} or to that of \cite{DaSa},\cite{Mi1},\cite{Mi2}.
By a simple transformation, this definition fits the rather general framework of 
$\bar\Gamma$-convergence, which is defined in Chap.~16 of \cite{Da}; 
that monograph however does not encompass the following Theorem~\ref{teo.comp}.  
\hfill$\Box$  
\end{remark}

Although we defined this convergence for a generic operator $\psi: L^p(0,T;X)\to L^1(0,T)$,
here we are concerned with superposition operators of the form  
\begin{equation}\label{eq.evol.super}
\begin{split}
&\psi_w(t) = \varphi(t,w(t))  
\qquad\forall w\in L^p_\mu(0,T;X),\hbox{ for a.e.\ }t\in {}]0,T[,
\\
&\varphi: {}]0,T[{} \times X\to \erre^+ \hbox{ \ being a {\it normal function.\/} }
\end{split}
\end{equation} 
By this we mean that $\varphi$ is globally measurable and $\varphi(t,\cdot)$ is lower semicontinuous for a.e.\ $t\in {}]0,T[$.  
  
\begin{theorem}[\cite{Vi17}]\label{teo.comp} 
(Compactness of evolutionary $\Gamma$-convergence of weak type)
Let $X$ be a real separable and reflexive Banach space, $p\in [1,+\infty[$, $T>0$,
and $\{\varphi_n\}$ be a sequence of normal functions ${}]0,T[{} \times X\to \erre^+$.
Let us assume that this sequence is equi-coercive and equi-bounded, in the sense that  
\begin{equation}\label{eq.evol.equibc} 
\begin{split}
&\exists C_1,C_2,C_3 >0: \forall n,\hbox{for a.e.\ }t\in {}]0,T[, \forall w\in X, 
\\
&C_1 \|w\|_X^p\le \varphi_n(t,w) \le C_2\|w\|_X^p +C_3, 
\end{split} 
\end{equation}  
and that
\begin{equation}\label{eq.evol.nul}
\varphi_n(t,0) =0 \qquad\hbox{ for a.e.\ }t\in {}]0,T[,\forall n.
\end{equation}  
\indent
Let $\mu$ be as in \eqref{eq.comstab.mu}. Let
$\tau$ be a topology on $L^p_\mu(0,T;X)$ that either coincides or is finer 
than the weak topology, and such that 
\begin{equation}\label{eq.evol.gamcom} 
\begin{split}
&\hbox{ for any sequence $\{F_n\}$ of functionals }
L^p_\mu(0,T;X)\to \erre^+\cup \{+\infty\},
\\
&\hbox{ if \ }\sup_{n\in\enne} \big\{\|w\|_{L^p_\mu(0,T;X)}: 
w\in L^p_\mu(0,T;X), F_n(w) \le C \big\}<+\infty,
\\
&\hbox{ then \ $\{F_n\}$ has a sequentially $\Gamma\tau$-convergent subsequence. }  
\end{split}
\end{equation}
\indent
Then there exists a normal function $\varphi: {}]0,T[{} \times X\to \erre^+$ such that
$\varphi(\cdot,0) =0$ a.e.\ in $]0,T[$, and such that, 
defining the operators $\psi,\psi_n: L^2_\mu(0,T;X)\to L^1_\mu(0,T)$ for any $n$ 
as in \eqref{eq.evol.super}, possibly extracting a subsequence,
\eqref{eq.evol.defgamma.1+} is fulfilled. 
\\
\indent
Moreover, if $\varphi_n$ does not depend on $t$ for any $n$, 
then the same holds for $\varphi$.
\end{theorem} 

\begin{remarks}\rm
(i) The integral representation of the $\Gamma$-limit is the main issue of this theorem;
this excluded the onset of long memory in the asymptotic functional.

(ii) The hypothesis  \eqref{eq.evol.gamcom} is crucial.
For instance, it is fulfilled by the weak topology, see e.g.\ chapter 8 of \cite{Da}.
In the case of interest of the present paper it is fulfilled, too, as we saw in Section~\ref{sec.comstab}. 
\hfill$\Box$ 
\end{remarks}

\bigskip 
\centerline{\bf Acknowledgment}
\medskip

The author is a member of GNAMPA of INdAM. 

This research was partially supported by a MIUR-PRIN 2015 grant for the project 
``Calcolo delle Variazioni" (Protocollo 2015PA5MP7-004).

\baselineskip=12truept
\bigskip

Author's address: 
\medskip

Augusto Visintin \par 
Universit\`a degli Studi di Trento \par 
Dipartimento di Matematica \par 
via Sommarive 14, \ 38050 Povo (Trento) - Italia \par 
Tel   +39-0461-281635 (office), +39-0461-281508 (secretary) \par 
Fax      +39-0461-281624 \par 
Email:   augusto.visintin@unitn.it 


\baselineskip=10.truept
\begin{thebibliography}{00} 

\bibitem{At}% 
{\sc H. Attouch:} 
{\it Variational Convergence for Functions and Operators.\/}
Pitman, Boston 1984

\bibitem{Au}% 
{\sc J.-P. Aubin:}
{\it Un th\'eor\`eme de compacit\'e.\/}
C. R. Acad. Sci. Paris {\bf 256} (1963) 5042--5044

\bibitem{AuEk}% 
{\sc J.-P. Aubin, I. Ekeland:}
{\it Applied Nonlinear Analysis.\/}
Wiley and Sons, New York 1984

\bibitem{AuFr}% 
{\sc J.-P. Aubin, H. Frankowska:}
{\it Set-valued analysis.\/} 
Birkh\"auser, Boston 1990 

\bibitem{Auc}%
{\sc G. Auchmuty:} 
{\it Saddle-points and existence-uniqueness for evolution equations.\/}
Differential Integral Equations {\bf 6} (1993) 1161--117  

\bibitem{BaBoWa}%
{\sc H.H. Bauschke, L.M. Borwein, X. Wang:} 
{\it Fitzpatrick functions and continuous linear monotone operators. \/}
SIAM J. Optim. {\bf 18} (2007) 789--809

\bibitem{BaWa}%
{\sc H.H. Bauschke, X. Wang:} 
{\it The kernel average for two convex functions and its applications to the extension and representation of monotone operators. \/}
Trans. Amer. Math. Soc. {\bf 361} (2009) 5947--5965

\bibitem{Bar}%
{\sc V. Barbu:} 
{\it Nonlinear Differential Equations of Monotone Types in Banach Spaces.\/} 
Springer, Berlin 2010

\bibitem{Bor}%
{\sc J.M. Borwein:} 
{\it Maximal monotonicity via convex analysis.\/}
J. Convex Anal. {\bf 13} (2006) 561--586

\bibitem{Bra1}%
{\sc A. Braides:} 
{\it $\Gamma$-Convergence for Beginners.\/}
Oxford University Press, Oxford 2002

\bibitem{Bra2}%
{\sc A. Braides:} 
{\it A Handbook of $\Gamma$-Convergence.\/}
In: Handbook of Partial Differential Equations.Stationary Partial Differential Equations, 
vol.\ 3 (M. Chipot, P. Quittner, Eds.) Elsevier, Amsterdam 2006, pp.\ 101Ð213

\bibitem{Br0} 
{\sc H. Brezis:} 
{\it \'Equations et in\'equations non lin\'eaires dans les espaces vectoriels en dualit\'e.\/}
Ann. Inst. Fourier (Grenoble) {\bf 18} (1968) 115--175

\bibitem{Br1}%
{\sc H. Brezis:} 
{\it Op\'erateurs Maximaux Monotones et Semi-Groupes de Contractions dans les
Espaces de Hilbert.\/}
North-Holland, Amsterdam 1973

\bibitem{BrEk}%
{\sc H. Brezis, I. Ekeland:} 
{\it Un principe variationnel associ\'e \`a certaines \'equations parabo\-liques. 
I. Le cas ind\'ependant du temps\/} 
and {\it II. Le cas d\'ependant du temps.\/}
C. R. Acad. Sci. Paris S\'er. A-B {\bf 282} (1976) 971--974,  and ibid. 1197--1198 

\bibitem{Bro1}%
{\sc F. Browder:}
{\it Mapping theorems for noncompact nonlinear operators in Banach spaces.\/}
Proc. Nat. Acad. Sci. U.S.A. {\bf 54} (1965) 337--342

\bibitem{Bro2}%
{\sc F. Browder:}
{\it Nonlinear operators and nonlinear equations of evolution in Banach spaces.\/}
Amer. Math. Soc., Providence, R. I. (1976) 1--308

\bibitem{BrHe} 
{\sc F. Browder and P. Hess:}
{\it Nonlinear mappings of monotone type in Banach spaces.\/}
J. Functional Analysis {\bf 11} (1972) 251--294
 
\bibitem{BuSv02}%
{\sc R.S. Burachik, B.F. Svaiter:}
{\it Maximal monotone operators, convex functions, and a special family of enlargements.\/} 
Set-Valued Analysis {\bf 10} (2002) 297--316  
 
\bibitem{BuSv03}%
{\sc R.S. Burachik, B.F. Svaiter:}
{\it Maximal monotonicity, conjugation and the duality product.\/}
Proc. Amer. Math. Soc. {\bf 131} (2003) 2379--2383  

\bibitem{CoVi}%
{\sc P. Colli, A. Visintin:} 
{\it On a class of doubly nonlinear evolution problems.\/}
Communications in P.D.E.s {\bf 15} (1990) 737--756
 
\bibitem{Da}%
{\sc G. Dal Maso:}
{\it An Introduction to $\Gamma$-Convergence.\/}
Birkh\"auser, Boston 1993  

\bibitem{DaSa}%
{\sc S. Daneri, G. Savar\`e:}
{\it Lecture notes on gradient flows and optimal transport.\/}
arXiv:1009.3737v1, 2010

\bibitem{DeFr}%
{\sc E. De Giorgi, T. Franzoni:}
{\it Su un tipo di convergenza variazionale.\/}
Atti Accad. Naz. Lincei Rend. Cl. Sci. Fis. Mat. Natur. (8) {\bf 58} (1975) 842--850 

\bibitem{DiSh}%
{\sc E. DiBenedetto, R.E. Showalter:} 
{\it Implicit degenerate evolution equations and applications.\/} 
S.I.A.M. J. Math. Anal. {\bf 12} (1981) 731--751 

\bibitem{DuLi}%
{\sc G. Duvaut, J.L. Lions:} 
{\it Les In\'equations en M\'ecanique et en Physique.\/} 
Dunod, Paris 1972

\bibitem{EkTe}%
{\sc I. Ekeland, R. Temam:}
{\it Analyse Convexe et Probl\`emes Variationnelles.\/}
Dunod Gau\-thier-Villars, Paris 1974

\bibitem{Fe}%
{\sc W. Fenchel:}
{\it Convex Cones, Sets, and Functions.\/}
Princeton Univ., 1953

\bibitem{Fi}%
{\sc S. Fitzpatrick:}
{\it Representing monotone operators by convex functions.\/}
Workshop/Mini\-con\-ference on Functional Analysis and Optimization (Canberra, 1988),
59--65,  Proc. Centre Math. Anal. Austral. Nat. Univ., 20, 
Austral. Nat. Univ., Canberra, 1988 

\bibitem{HuPa}
{\sc S. Hu, N. S. Papageorgiou:}
{\it Handbook of Multivalued Analysis.\/} 
Vol.\ I, Kluwer, Dordrecht 1979  

\bibitem{Ke}%
{\sc N. Kenmochi:} 
{\it Monotonicity and compactness methods for nonlinear variational inequalities.\/}
In: Handbook of differential equations: stationary partial differential equations. 
Vol. IV, 203--298, Handb. Differ. Equ., Elsevier/North-Holland, Amsterdam, 2007 

\bibitem{Kr}%
{\sc N.V. Krylov:}
{\it Some properties of monotone mappings.\/}
Litovsk. Mat. Sb. {\bf 22} (1982) 80--87 

\bibitem{Li}%
{\sc J.L. Lions:}
{\it Quelques M\'ethodes de R\'esolution des Probl\`emes aux Limites non Lin\'eaires.\/}
Dunod, Paris 1969

\bibitem{MaRa}%
{\sc J. M\'alek, K.R. Rajagopal:}
{\it Mathematical issues concerning the Navier-Stokes equations and some 
of its generalizations.\/}
In: Evolutionary equations. Handb. Differ. Equ., 
Elsevier/North-Holland, Amsterdam (2005). Vol. II, pp.\ 371--459

\bibitem{MaSv05}%
{\sc J.-E. Martinez-Legaz, B.F. Svaiter:}
{\it Monotone operators representable by l.s.c.\ convex functions.\/}
Set-Valued Anal. {\bf 13} (2005) 21--46 

\bibitem{MaSv08}%
{\sc J.-E. Martinez-Legaz, B.F. Svaiter:}
{\it Minimal convex functions bounded below by the duality product.\/} 
Proc. Amer. Math. Soc. {\bf 136} (2008) 873--878  

\bibitem{MaTh01}%
{\sc J.-E. Martinez-Legaz, M. Th\'era:}
{\it A convex representation of maximal monotone operators.\/}  
J. Nonlinear Convex Anal. {\bf 2} (2001) 243--247 

\bibitem{Mi1}%
{\sc A. Mielke:} 
{\it Lecture notes on evolutionary $\Gamma$-convergence for gradient systems.\/}
In: Macroscopic and Large Scale Phenomena: Coarse Graining, Mean Field Limits and Ergodicity (A. Muntean, J. Rademacher, A. Zagaris, eds.).
Lecture Notes in Appl. Math. Mech.. Springer 2016, pp.\ 187--249 

\bibitem{Mi2}%
{\sc A. Mielke:}
{\it Deriving amplitude equations via evolutionary $\Gamma$-convergence.\/}
Discrete Contin. Dyn. Syst. Ser. A, {\bf 35} (2015) 2679--2700  

\bibitem{Mo}%
{\sc U. Mosco:}  
{\it Convergence of convex sets and of solutions of variational inequalities.\/}
Advances in Math. {\bf 3} (1969) 510--585 

\bibitem{Na}%
{\sc B. Nayroles:}
{\it Deux th\'eor\`emes de minimum pour certains syst\`emes dissipatifs.\/}
C. R. Acad. Sci. Paris S\'er. A-B {\bf 282} (1976) A1035--A1038 

\bibitem{PaSb}%
{\sc D. Pascali, S. Sburlan:}
{\it Nonlinear mappings of monotone type.\/}
Martinus Nijhoff, The Hague; Sijthoff and Noordhoff, Alphen aan den Rijn, 1978

\bibitem{Pe04}%
{\sc J.-P. Penot:}
{\it A representation of maximal monotone operators by closed convex functions and its   
impact on calculus rules.\/} 
C. R. Math. Acad. Sci. Paris, Ser.\ I {\bf 338} (2004) 853--858

\bibitem{Pe04'}% 
{\sc J.-P. Penot:}
{\it The relevance of convex analysis for the study of monotonicity.\/}
Nonlinear Anal. {\bf 58} (2004) 855--871

\bibitem{Ri1}%ÊÊ
{\sc H. Rios:}
{\it \'Etude de la question d'existence pour certains probl\`emes d'\'evolution par
minimisation d'une fonctionnelle convexe.\/}
C. R. Acad. Sci. Paris S\'er. A-B {\bf 283} (1976) A83--A86Ê 

\bibitem{RoRoSt}%
{\sc T. Roche, R. Rossi, U. Stefanelli:}
{\it Stability results for doubly nonlinear differential inclusions by variational convergence.\/}
SIAM J. Control Optim. {\bf 52} (2014) 1071--1107 

\bibitem{Rou2}% 
{\sc T. Roub\'\i\v cek:} 
{\it Durect method for parabolic problems.\/}
Adv. Math. Sci. Appl. {\bf 10} (2000) 57--65

\bibitem{SaSe}%
{\sc E. Sandier, S. Serfaty:}
{\it Gamma-convergence of gradient flows and applications to Ginzburg-Landau.\/} 
Comm. Pure Appl. Math. {\bf 55} (2002) 537--581  

\bibitem{St1}%
{\sc U. Stefanelli:}
{\it The Brezis-Ekeland principle for doubly nonlinear equations.\/} 
SIAM J. Control Optim. {\bf 47} (2008) 1615--1642 

\bibitem{Sv03}%
{\sc B.F. Svaiter:}
{\it Fixed points in the family of convex representations of a maximal monotone operator.\/} 
Proc. Amer. Math. Soc. {\bf 131} (2003) 3851--3859 

\bibitem{Te}%
{\sc R. Temam:}
{\it Navier-Stokes equations. Theory and numerical analysis.\/} 
North-Holland, Amsterdam-New York 1979
 
\bibitem{Vi96}%
{\sc A. Visintin:} 
{\it Models of Phase Transitions.\/}
Birkh\"auser, Boston 1996  

\bibitem{ViAMSA}%
{\sc A. Visintin:}
{\it Extension of the Brezis-Ekeland-Nayroles principle to monotone operators.\/}
Adv. Math. Sci. Appl. {\bf 18} (2008) 633--650 

\bibitem{ViCalVar}%
{\sc A. Visintin:}
{\it Variational formulation and structural stability of monotone equations.\/}  
Calc. Var. Partial Differential Equations {\bf 47} (2013) 273--317 

\bibitem{Vi14}%
{\sc A. Visintin:}
{\it An extension of the Fitzpatrick theory.\/}  
Commun.\ Pure Appl. Anal. {\bf 13} (2014) 2039--2058

\bibitem{Vi15}%
{\sc A. Visintin:}
{\it Weak structural stability of pseudo-monotone equations.\/}
Discrete Contin. Dyn. Syst. Ser. A  35 (2015) 2763--2796 

\bibitem{ViGil}%
{\sc A. Visintin:}
{\it On the structural properties of nonlinear flows.\/}  
In: ``Solvability, Regularity, Optimal Control of Boundary Value Problems for PDEs''
Springer INdAM Series
(in press)  

\bibitem{Vi17}%
{\sc A. Visintin:}
{\it Evolutionary $\Gamma$-convergence of weak type.\/}  
(arXiv) 


\end{thebibliography}
\end{document}